\documentclass[a4paper]{article}

\usepackage{amsfonts,amssymb,amsthm}

\makeatletter
\newif\iflabel\labelfalse \let\w@label=\label
\def\label{\global\labeltrue\w@label}
\let\w@eqn=\equation
\def\equation{\global\labelfalse\w@eqn}
\def\endequation{\iflabel\@eeq\else\@eeqw\fi$$\global\@ignoretrue}
\def\@eeq{\eqno \@eqnnum}
\def\@eeqw{\addtocounter{equation}{-1}}
\newskip\bw@\newskip\bws@
\def\dc@@pq{\global \bw@=\belowdisplayskip \global\bws@=\belowdisplayshortskip
   \global\belowdisplayshortskip=3pt plus -3pt
   \global\belowdisplayskip=\belowdisplayshortskip
   \hskip 500pt minus 500pt\relax}
\def\dc@pq{\dc@@pq$$}
\def\fc@pq{$$\hskip 500pt minus 500pt\global
   \belowdisplayshortskip=\bws@\global\belowdisplayskip=\bw@}
\def\coupeq{\dc@pq\fc@pq}
\def\coupepas{\par\noindent\begin{minipage}{\textwidth}}
\def\jusquela{\end{minipage}}
\makeatother

\newcommand{\dr}{\partial}

\newcommand{\C}{{\mathbb C}}
\newcommand{\N}{{\mathbb N}}
\newcommand{\R}{{\mathbb R}}
\newcommand{\Z}{{\mathbb Z}}
\newcommand{\II}{I\hspace{-0.1cm}I}

\newcommand{\argsinh}{\mbox{argsinh}}

\newcommand{\hess}{\mbox{Hess}}

\newcommand{\rhob}{\overline{\rho}}

\newcommand{\St}{\tilde{S}}

\newtheorem{prop}{Proposition}[section]
\newtheorem{lemma}[prop]{Lemma}
\newtheorem{sublemma}[prop]{Sub-lemma}

\newtheorem{thm}[prop]{Theorem}
\newtheorem{cor}[prop]{Corollary}
\newtheorem{remark}[prop]{Remark}

\newtheorem{df}[prop]{Definition}
\newtheorem{pty}[prop]{Property}
\newtheorem{question}[prop]{Question}

\newcommand{\pg}{\paragraph}

\newenvironment{thn}[1]{\vskip 0.2cm \noindent{\bf Theorem #1.}\it}{\rm
\vspace{0.2cm}} 
 
\newenvironment{crn}[1]{\vskip 0.2cm \noindent{\bf Corollary #1.}\it}{\rm
\vspace{0.2cm}} 
\newenvironment{lmn}[1]{\vskip 0.2cm \noindent{\bf Lemma #1.}\it}{\rm
\vspace{0.2cm}} 
\newenvironment{qn}[1]{\vskip 0.2cm \noindent{\bf Question #1.}\it}{\rm
\vspace{0.2cm}} 

\newcommand{\btm}{\begin{thm}}
\newcommand{\etm}{\end{thm}}
\newcommand{\bpt}{\begin{pty}}
\newcommand{\ept}{\end{pty}}
\newcommand{\blm}{\begin{lemma}}
\newcommand{\elm}{\end{lemma}}
\newcommand{\bsl}{\begin{sublemma}}
\newcommand{\esl}{\end{sublemma}}
\newcommand{\bcr}{\begin{cor}}
\newcommand{\ecr}{\end{cor}}
\newcommand{\bdf}{\begin{df}}
\newcommand{\edf}{\end{df}}
\newcommand{\bprop}{\begin{prop}}
\newcommand{\eprop}{\end{prop}}
\newcommand{\bas}{\begin{asser}}
\newcommand{\eas}{\end{asser}}
\newcommand{\beq}{\begin{equation}}
\newcommand{\eeq}{\end{equation}}
\newcommand{\bpv}{\begin{proof}}
\newcommand{\epv}{\end{proof}}
\newcommand{\bpvs}{\begin{sketch}}
\newcommand{\epvs}{\end{sketch}}
\newcommand{\bit}{\begin{itemize}}
\newcommand{\eit}{\end{itemize}}
\newcommand{\bpn}{\begin{pfn}}
\newcommand{\epn}{\end{pfn}}
\newcommand{\btn}{\begin{thn}}
\newcommand{\etn}{\end{thn}}
\newcommand{\bcn}{\begin{crn}}
\newcommand{\ecn}{\end{crn}}
\newcommand{\bqn}{\begin{qn}}
\newcommand{\eqn}{\end{qn}}
\newcommand{\bln}{\begin{lmn}}
\newcommand{\eln}{\end{lmn}}
\newcommand{\brk}{\begin{remark}}
\newcommand{\erk}{\end{remark}}
\newcommand{\bq}{\begin{question}}
\newcommand{\eq}{\end{question}}

\newenvironment{pfn}[1]{\vskip 0.2cm \noindent{\it Proof #1.}}{$\square$
\vspace{0.2cm}}

\newcommand{\cA}{{\mathcal A}}

\newcommand{\cP}{{\mathcal P}}

\newcommand{\cT}{{\mathcal T}}

\newcommand{\Isom}{\mathrm{Isom}}

\newcommand{\ud}{\dot{u}}
\newcommand{\vd}{\dot{v}}

\newcommand{\alphad}{\dot{\alpha}}
\newcommand{\phid}{\dot{\phi}}
\newcommand{\rhod}{\dot{\rho}}
\newcommand{\thetad}{\dot{\theta}}

\catcode`\@=11\relax
\newwrite\@unused
\def\typeout#1{{\let\protect\string\immediate\write\@unused{#1}}}
\def\@nnil{\@nil}
\def\@empty{}
\def\@psdonoop#1\@@#2#3{}
\def\@psdo#1:=#2\do#3{\edef\@psdotmp{#2}\ifx\@psdotmp\@empty \else
    \expandafter\@psdoloop#2,\@nil,\@nil\@@#1{#3}\fi}
\def\@psdoloop#1,#2,#3\@@#4#5{\def#4{#1}\ifx #4\@nnil \else
       #5\def#4{#2}\ifx #4\@nnil \else#5\@ipsdoloop #3\@@#4{#5}\fi\fi}
\def\@ipsdoloop#1,#2\@@#3#4{\def#3{#1}\ifx #3\@nnil 
       \let\@nextwhile=\@psdonoop \else
      #4\relax\let\@nextwhile=\@ipsdoloop\fi\@nextwhile#2\@@#3{#4}}
\def\@tpsdo#1:=#2\do#3{\xdef\@psdotmp{#2}\ifx\@psdotmp\@empty \else
    \@tpsdoloop#2\@nil\@nil\@@#1{#3}\fi}
\def\@tpsdoloop#1#2\@@#3#4{\def#3{#1}\ifx #3\@nnil 
       \let\@nextwhile=\@psdonoop \else
      #4\relax\let\@nextwhile=\@tpsdoloop\fi\@nextwhile#2\@@#3{#4}}
\def\psdraft{
        \def\@psdraft{0}
}
\def\psfull{
        \def\@psdraft{100}
}
\psfull
\newif\if@prologfile
\newif\if@postlogfile
\newif\if@noisy
\def\pssilent{
        \@noisyfalse
}
\def\psnoisy{
        \@noisytrue
}
\psnoisy
\newif\if@bbllx
\newif\if@bblly
\newif\if@bburx
\newif\if@bbury
\newif\if@height
\newif\if@width
\newif\if@rheight
\newif\if@rwidth
\newif\if@clip
\newif\if@verbose
\def\@p@@sclip#1{\@cliptrue}
\def\@p@@sfile#1{%\typeout{file is #1}
                   \def\@p@sfile{#1}
}
\def\@p@@sfigure#1{\def\@p@sfile{#1}}
\def\@p@@sbbllx#1{
                \@bbllxtrue
                \dimen100=#1
                \edef\@p@sbbllx{\number\dimen100}
}
\def\@p@@sbblly#1{
                \@bbllytrue
                \dimen100=#1
                \edef\@p@sbblly{\number\dimen100}
}
\def\@p@@sbburx#1{
                \@bburxtrue
                \dimen100=#1
                \edef\@p@sbburx{\number\dimen100}
}
\def\@p@@sbbury#1{
                \@bburytrue
                \dimen100=#1
                \edef\@p@sbbury{\number\dimen100}
}
\def\@p@@sheight#1{
                \@heighttrue
                \dimen100=#1
                \edef\@p@sheight{\number\dimen100}
}
\def\@p@@swidth#1{
                \@widthtrue
                \dimen100=#1
                \edef\@p@swidth{\number\dimen100}
}
\def\@p@@srheight#1{
                \@rheighttrue
                \dimen100=#1
                \edef\@p@srheight{\number\dimen100}
}
\def\@p@@srwidth#1{
                \@rwidthtrue
                \dimen100=#1
                \edef\@p@srwidth{\number\dimen100}
}
\def\@p@@ssilent#1{ 
                \@verbosefalse
}
\def\@p@@sprolog#1{\@prologfiletrue\def\@prologfileval{#1}}
\def\@p@@spostlog#1{\@postlogfiletrue\def\@postlogfileval{#1}}
\def\@cs@name#1{\csname #1\endcsname}
\def\@setparms#1=#2,{\@cs@name{@p@@s#1}{#2}}
\def\ps@init@parms{
                \@bbllxfalse \@bbllyfalse
                \@bburxfalse \@bburyfalse
                \@heightfalse \@widthfalse
                \@rheightfalse \@rwidthfalse
                \def\@p@sbbllx{}\def\@p@sbblly{}
                \def\@p@sbburx{}\def\@p@sbbury{}
                \def\@p@sheight{}\def\@p@swidth{}
                \def\@p@srheight{}\def\@p@srwidth{}
                \def\@p@sfile{}
                \def\@p@scost{10}
                \def\@sc{}
                \@prologfilefalse
                \@postlogfilefalse
                \@clipfalse
                \if@noisy
                        \@verbosetrue
                \else
                        \@verbosefalse
                \fi
}
\def\parse@ps@parms#1{
                \@psdo\@psfiga:=#1\do
                   {\expandafter\@setparms\@psfiga,}}
\newif\ifno@bb
\newif\ifnot@eof
\newread\ps@stream
\def\bb@missing{
        \if@verbose{
                \typeout{psfig: searching \@p@sfile \space  for bounding box}
        }\fi
        \openin\ps@stream=\@p@sfile
        \no@bbtrue
        \not@eoftrue
        \catcode`\%=12
        \loop
                \read\ps@stream to \line@in
                \global\toks200=\expandafter{\line@in}
                \ifeof\ps@stream \not@eoffalse \fi
                \@bbtest{\toks200}
                \if@bbmatch\not@eoffalse\expandafter\bb@cull\the\toks200\fi
        \ifnot@eof \repeat
        \catcode`\%=14
}       
\catcode`\%=12
\newif\if@bbmatch
\def\@bbtest#1{\expandafter\@a@\the#1%%BoundingBox:\@bbtest\@a@}
\long\def\@a@#1%%BoundingBox:#2#3\@a@{\ifx\@bbtest#2\@bbmatchfalse\else\@bbmatchtrue\fi}
\long\def\bb@cull#1 #2 #3 #4 #5 {
        \dimen100=#2 bp\edef\@p@sbbllx{\number\dimen100}
        \dimen100=#3 bp\edef\@p@sbblly{\number\dimen100}
        \dimen100=#4 bp\edef\@p@sbburx{\number\dimen100}
        \dimen100=#5 bp\edef\@p@sbbury{\number\dimen100}
        \no@bbfalse
}
\catcode`\%=14
\def\compute@bb{
                \no@bbfalse
                \if@bbllx \else \no@bbtrue \fi
                \if@bblly \else \no@bbtrue \fi
                \if@bburx \else \no@bbtrue \fi
                \if@bbury \else \no@bbtrue \fi
                \ifno@bb \bb@missing \fi
                \ifno@bb \typeout{FATAL ERROR: no bb supplied or found}
                        \no-bb-error
                \fi
                \count203=\@p@sbburx
                \count204=\@p@sbbury
                \advance\count203 by -\@p@sbbllx
                \advance\count204 by -\@p@sbblly
                \edef\@bbw{\number\count203}
                \edef\@bbh{\number\count204}
}
\def\in@hundreds#1#2#3{\count240=#2 \count241=#3
                     \count100=\count240        % 100 is first digit #2/#3
                     \divide\count100 by \count241
                     \count101=\count100
                     \multiply\count101 by \count241
                     \advance\count240 by -\count101
                     \multiply\count240 by 10
                     \count101=\count240        %101 is second digit of #2/#3
                     \divide\count101 by \count241
                     \count102=\count101
                     \multiply\count102 by \count241
                     \advance\count240 by -\count102
                     \multiply\count240 by 10
                     \count102=\count240        % 102 is the third digit
                     \divide\count102 by \count241
                     \count200=#1\count205=0
                     \count201=\count200
                        \multiply\count201 by \count100
                        \advance\count205 by \count201
                     \count201=\count200
                        \divide\count201 by 10
                        \multiply\count201 by \count101
                        \advance\count205 by \count201
                     \count201=\count200
                        \divide\count201 by 100
                        \multiply\count201 by \count102
                        \advance\count205 by \count201
                     \edef\@result{\number\count205}
}
\def\compute@wfromh{
                \in@hundreds{\@p@sheight}{\@bbw}{\@bbh}
                \edef\@p@swidth{\@result}
}
\def\compute@hfromw{
                \in@hundreds{\@p@swidth}{\@bbh}{\@bbw}
                \edef\@p@sheight{\@result}
}
\def\compute@handw{
                \if@height 
                        \if@width
                        \else
                                \compute@wfromh
                        \fi
                \else 
                        \if@width
                                \compute@hfromw
                        \else
                                \edef\@p@sheight{\@bbh}
                                \edef\@p@swidth{\@bbw}
                        \fi
                \fi
}
\def\compute@resv{
                \if@rheight \else \edef\@p@srheight{\@p@sheight} \fi
                \if@rwidth \else \edef\@p@srwidth{\@p@swidth} \fi
}
\def\compute@sizes{
        \compute@bb
        \compute@handw
        \compute@resv
}
\def\psfig#1{\vbox {
        \ps@init@parms
        \parse@ps@parms{#1}
        \compute@sizes
        \ifnum\@p@scost<\@psdraft{
                \if@verbose{
                        \typeout{psfig: including \@p@sfile \space }
                }\fi
                \special{ps::[begin]    \@p@swidth \space \@p@sheight \space
                                \@p@sbbllx \space \@p@sbblly \space
                                \@p@sbburx \space \@p@sbbury \space
                                startTexFig \space }
                \if@clip{
                        \if@verbose{
                                \typeout{(clip)}
                        }\fi
                        \special{ps:: doclip \space }
                }\fi
                \if@prologfile
                    \special{ps: plotfile \@prologfileval \space } \fi
                \special{ps: plotfile \@p@sfile \space }
                \if@postlogfile
                    \special{ps: plotfile \@postlogfileval \space } \fi
                \special{ps::[end] endTexFig \space }
                \vbox to \@p@srheight true sp{
                        \hbox to \@p@srwidth true sp{
                                \hss
                        }
                \vss
                }
        }\else{
                \vbox to \@p@srheight true sp{
                \vss
                        \hbox to \@p@srwidth true sp{
                                \hss
                                \if@verbose{
                                        \@p@sfile
                                }\fi
                                \hss
                        }
                \vss
                }
        }\fi
}}
\catcode`\@=12\relax

\begin{document}

\title{Small deformations of polygons}

\author{Jean-Marc Schlenker\thanks{
Laboratoire Emile Picard, UMR CNRS 5580,
UFR MIG, Universit{\'e} Paul Sabatier,
118 route de Narbonne,
31062 Toulouse Cedex 4,
France.
\texttt{schlenker@picard.ups-tlse.fr; http://picard.ups-tlse.fr/\~{
}schlenker}. }}

\date{October 2004; revised, December 2004 (v2)}

\maketitle

\begin{abstract}
We describe the first-order variations of the angles of Euclidean,
spherical or 
hyperbolic polygons under infinitesimal deformations such that the lengths of
the edges do not change. Using 
this description, we introduce a vector-valued 
quadratic invariant $b$ on the space of
those isometric deformations which, for 
convex polygons, has a remarkable positivity
property.

We give two geometric applications. The first is an isoperimetric statement
for hyperbolic polygons: among the convex hyperbolic 
polygons with given edge lengths, there is a unique polygon
with vertices on a circle, a horocycle, or on one connected component of the
space of points at constant
distance from a geodesic, and it has maximal area. 
The second application is a new proof of the infinitesimal rigidity of convex
polyhedra in the Euclidean space, and a new rigidity result for
polyhedral surfaces in the Minkowski space. 

Finally we indicate how the invariant $b$ can 
be used to define natural metrics on the space of
convex spherical (or hyperbolic) polygons with fixed edge lengths. Those
metrics are related to known (and interesting) metrics on the space of
convex polygons with given angles in the plane. 
\bigskip

\begin{center} {\bf R{\'e}sum{\'e}} \end{center}
On d\'ecrit les d\'eformations infinit\'esimales des angles 
d'un polygone euclidien, sph\'erique
ou hyperbolique sous les d\'eformations infinit\'esimales 
qui pr\'eservent les longueurs des ar\^etes. On en d\'eduit la
d\'efinition d'un invariant quadratique \`a valeurs vectorielles 
$b$ sur l'espace de ces d\'eformations isom\'etriques qui, 
pour les polygones convexes, a une propri\'et\'e
remarquable de positivit\'e.

On donne deux applications g\'eom\'etriques. La premi\`ere est un \'enonc\'e
isoperim\'etrique pour les polygones
hyperboliques: parmi les polygones hyperboliques convexes dont les longueurs
des ar\^etes sont donn\'ees, il existe 
un unique \'el\'ement dont les sommets sont sur un
cercle, un horocycle, ou dans une composante connexe de l'ensemble des points 
\`a distance constante d'une droite, et son aire est
maximale. La seconde application est une nouvelle preuve de la rigidit\'e
infinit\'esimale des poly\`edres euclidiens, et un nouveau r\'esultat de
rigidit\'e pour les surfaces poly\`edrales dans l'espace de Minkowski.

Finalement on indique comment l'invariant $b$ peut \^etre utilis\'e pour
d\'efinir des m\'etriques naturelles sur l'espace des polygones convexes
sph\'eriques (ou hyperboliques) dont les longueurs des cot\'es sont
fix\'es. Ces m\'etriques sont reli\'ees \`a des m\'etriques connues, et
int\'ressantes, sur les espaces de polygones euclidiens convexes dont les
angles sont fix\'es.

\end{abstract}

\tableofcontents

\section{Introduction}

\pg{Moduli spaces of polygons.} 

In this paper, a polygon with $n$ vertices in
the Euclidean plane $E^2$ (resp. the sphere $S^2$, the hyperbolic plane $H^2$
or the de Sitter plane $S^2_1$)
is a finite sequence of points $v_1, v_2, \cdots v_n$ in $E^2$ (resp. $S^2$,
$H^2$, $S^2_1$), with $v_0:=v_n$, and with $v_i\neq v_{i+1}$ for all $i\in
\{1,\cdots, n\}$. We call $\cP_{n,E}$ (resp.  $\cP_{n,S}$, 
$\cP_{n,H}$, $\cP_{n,dS}$) the moduli space of polygons in each of the spaces,
i.e. the space of polygons quotiented by the group of orientation-preserving
isometries of the space. In the sphere and the de Sitter plane, we only
consider polygons such that two consecutive vertices are never antipodal (in
the sphere, they are at distance less than $\pi$) although those cases could
be included at the cost of a little more care. 

The edge lengths of a polygon are the numbers 
$l_1, \cdots, l_n$, where $l_i$ is the distance between $v_i$ and
$v_{i+1}$. Given $l=(l_1, \cdots, l_n)$, $\cP_E(l)$ is the subset of
$\cP_{n,E}$ of Euclidean 
polygons with edge lengths equal to $l_1, \cdots, l_n$.
Similar notations will be used for spherical, hyperbolic or de Sitter
polygons. 

\pg{First-order deformations of polygons.}

Let $p=(v_1, \cdots, v_n)$ be a Euclidean polygon, and let $\alpha_1, \cdots,
\alpha_n$ be its angles. An isometric first-order deformation of $p$ is a set
of vectors $\vd_1, \cdots, \vd_n$, with $\vd_i\in T_{v_i}E^2$, such that, if
one deforms $p$ infinitesimally by moving each $v_i$ along $\vd_i$, the
lengths of the edges do not change (at first order). This is equivalent to
the fact that, for each $i\in \{ 0, 1,\cdots, n-1\}$, $\langle v_{i+1}-v_i,
\vd_{i+1}-\vd_i\rangle=0$. An infinitesimal first-order deformation of $p$ is
{\it trivial} if there exists a Killing field $\kappa$ on $\R^2$ such that,
for each $i$, $\vd_i=\kappa(v_i)$. We mostly consider the isometric
first-order deformations of $p$ up to the trivial deformations; for
a ``generic'' polygon $p$, they are
canonically associated to the elements of $T_{[p]}\cP_{E}(l)$, 
where $l=(l_1, \cdots, l_n)$ is the set of lengths of $p$ and $[p]$ is the
image of $p$ under the projection from the space of polygons in $\R^2$ of
given edge lengths to the quotient under the action of $\Isom(\R^2)$.

Clearly, the isometric first-order deformations of $p$ are determined by the
first-order variations $\alphad_i$ of its angles. However not all possible
variations are possible. We first give some constraints on the possible
first-order variations of the angles of polygons in the spaces we consider. 

\begin{prop} \label{pr:euclidien}
  Let $p=(v_1, \cdots, v_n)$ be a Euclidean polygon. Let $\alpha_1, \cdots,
  \alpha_n$ be its angles, and let $\alphad_1, \cdots, \alphad_n\in
  \R$. Suppose that there
  exists a first-order isometric deformation of $p$ such that the induced
  first-order variation of the angles of $p$ is $(\alphad_1, \cdots,
  \alphad_n)$. Then:
  \begin{itemize}
  \item $\sum_{i=1}^n \alphad_i=0$.
  \item $\sum_{i=1}^n \alphad_i v_i=0$.
  \end{itemize}
  Conversely, if those two conditions are satisfied and the vertices of $p$
  are not all on a line, then there exists an isometric 
  first-order deformation of $p$
  such that the induced 
  first-order variation of the $\alpha_i$ are the $\alphad_i$, and it is
  unique up to the addition of a trivial deformation.
\end{prop}

An analogous but even simpler statement holds for spherical polygons.

\begin{thn}{A$_S$}
  Let $p=(v_1, \cdots, v_n)$ be a spherical polygon. Let $\alpha_1, \cdots,
  \alpha_n$ be its angles, and let $\alphad_1, \cdots, \alphad_n\in \R$
  be a first-order variation of its angles induced by an isometric first-order
  deformation of $p$. Then:
$$ \sum_{i=1}^n \alphad_i v_i=0~, $$
where the $v_i$ are considered as points in $S^2\subset \R^3$. 

Conversely, if this
equation is satisfied by an $n$-uple $(\alphad_1, \cdots,\alphad_n)$ and
moreover the $v_i$ are not all on a great circle, then there exists an
isometric first-order deformation of $p$ such that the $\alphad_i$ are the
associated first-order variations of the $\alpha_i$.
\end{thn}

Each first-order isometric deformation of $p$ is uniquely determined, up to
the addition of a trivial deformation, by the
corresponding variation of its angles. Moreover, the space of isometric
first-order deformations of $p$ (up to the trivial ones)
has dimension at least $n-3$, since a
deformation of $p$ is determined by $2n$ parameters (the positions of the
vertices), with $n$ constraints (the lengths of the edges) and a group of
isometries of dimensions $3$. Therefore,
the space of first-order variations of the $\alpha_i$, induced by isometric
first-order deformations of $p$, always has dimension at least $n-3$. It
follows that, unless the $v_i$ are on a great circle, the space of isometric
first-order deformations of $p$ has dimension $n-3$. We will see at the
beginning of section 2 that this implies the 
fact, proved by Kapovich and Millson \cite{kapovich-millson-sphere},
that in this case the space of polygons with the same edge lengths as $p$ is
locally a smooth manifold. 

The same description holds in the hyperbolic plane. 

\begin{thn}{A$_H$}
  Let $p=(v_1, \cdots, v_n)$ be a hyperbolic polygon. Let $\alpha_1, \cdots,
  \alpha_n$ be its angles, and let $\alphad_1, \cdots, \alphad_n\in \R$
  be a first-order variation of its angles induced by an isometric first-order
  deformation of $p$. Then:
$$ \sum_{i=1}^n \alphad_i v_i=0~, $$
where the $v_i$ are considered as points in $H^2\subset \R^3_1$. 

Conversely, if this
equation is satisfied by an $n$-uple $(\alphad_1, \cdots,\alphad_n)$ and
moreover the $v_i$ are not all on a hyperbolic geodesic, then there exists an
isometric first-order deformation of $p$ such that the $\alphad_i$ are the
associated first-order variations of the $\alpha_i$.
\end{thn}

Here $\R^3_1$ is the Minkowski 3-dimensional space, i.e. it is $\R^3$ with the
bilinear form:
$$ \langle (x,y,z), (x',y',z')\rangle = xx' + yy' - zz'~. $$
$H^2$ has an isometric embedding in $\R^3_1$ as a quadric:
$$ H^2 = \{ (x,y,z)\in \R^3_1 ~ | ~ \langle (x,y,z),(x,y,z)\rangle = -1 ~
 \mbox{and} ~ z>0 \}~. $$

\pg{A quadratic form on first-order deformations.}

This description of the first-order deformations of the angles of polygons
opens the door to 
the definition of a quadratic invariant of first-order isometric
deformations. 

\bdf \label{df:bUS}
Let $p=(v_1, \cdots, v_n)$ be an Euclidean polygon, with edge lengths $l_1,
\cdots, l_n$. Let $l:=(l_1, \cdots, l_n)$, and let $U\in
T_{[p]}\cP_{E}(l)$, where $[p]$ is the projection of $p$ in $\cP_E(l)$. Let
$U'$ be a first-order deformation of $p$ projecting, under the quotient by
the trivial deformations, to $U$. We call: 
$$ b(U) := \sum_{i=1}^n d\alpha_i(U') dv_i(U')~, $$
where the $d\alpha_i(U')$ and the $dv_i(U')$ are the first-order variations of
the $\alpha_i$ and of the $v_i$, respectively, under $U'$. Then $b(U)$ does
not depend on the choice of $U'$. 
\edf

As a first-order deformation of $p$ defined up to the trivial deformations,
$U$ does not define uniquely the first-order variations of the $v_i$. However,
any choice $U'$ of a representative $(\vd_1, \cdots, \vd_n)$ will do. Indeed,
let $U''$ be another first-order deformation of $p$ corresponding to $U$; then
$U'-U''$ is a trivial deformation, so the corresponding variation of the $v_i$
is given by $\kappa(v_i)$, where $\kappa$ is a Killing field. 
But any Killing field $\kappa$ in $E^2$ is of the form: 
$$ \kappa(x) = e^{i\theta} x + \tau~, $$ 
for some $\theta \in \R$ and some vector $\tau\in E^2$ (this expression uses
the usual identification of $E^2$ with $\C$). It follows that: 
$$
\sum_{i=1}^n d\alpha_i(U')dv_i(U') - \sum_{i=1}^n d\alpha_i(U'') dv_i(U'') =
\sum_{i=1}^n d\alpha_i(U') \kappa(v_i) = \coupeq
= \sum_{i=1}^n d\alpha_i(U')
(e^{i\theta} v_i + \tau) 
= e^{i\theta} \left(\sum_{i=1}^n d\alpha_i(U')
  v_i\right) + \left(\sum_{i=1}^n d\alpha_i(U') \right) \tau~, 
$$
and both terms vanish by Proposition \ref{pr:euclidien}.

A very similar definition can be used in the sphere or the hyperbolic plane.

\bdf
Let $p=(v_1, \cdots, v_n)$ be a polygon in the sphere (resp. 
the hyperbolic plane) with edge lengths $l_1,
\cdots, l_n$. Let $l:=(l_1, \cdots, l_n)$, and let $U\in T_{[p]}\cP_{S}(l)$
(resp. $T_{[p]}\cP_{H}(l)$) where $[p]$ is the projection of $p$ in $\cP_S(l)$
(resp. $\cP_H(l)$). Let $U'$ be a first-order
deformation of $p$ projecting, under the quotient by the trivial
deformations, to $U$. We call: 
$$ b(U) := \sum_{i=1}^n d\alpha_i(U') dv_i(U')~, $$
where the $d\alpha_i(U')$ and the $dv_i(U')$ are the first-order variations of
the $\alpha_i$ and of the $v_i$, considered as points in $\R^3$
(resp. $\R^3_1$), respectively, under $U'$. Then $b(U)$ does not depend on the
choice of $U'$.
\edf

An argument very similar to the one given above shows that $b(U)$ is indeed
independent of the precise deformation $U'$ of $p$ which is chosen. 
In the sphere, any Killing field is of the form:
$$ \kappa(x) = Y\times x~, $$
for some vector $Y\in \R^3$. Therefore, if $U''$ is another first-order
deformation of $p$ projecting to $U$, then $U'-U''$ corresponds to a
Killing field $\kappa$. Then:
$$
\sum_{i=1}^n d\alpha_i(U')dv_i(U') - \sum_{i=1}^n d\alpha_i(U'')dv_i(U'')  
=  \sum_{i=1}^n d\alpha_i(U') \kappa(v_i) = \coupeq 
= \sum_{i=1}^n d\alpha_i(U') Y\times v_i 
= Y\times \left(\sum_{i=1}^n d\alpha_i(U') v_i \right) 
= 0~. 
$$
The same argument can be used in the hyperbolic context,
using the fact that the Killing fields of $\R^3_1$ which vanish at the origin
(in other terms the elements of the Lie algebra $so(2,1)$) are of the form:
$$ \kappa(x) = Y\boxtimes x~, $$
where $\boxtimes$ is the Minkowski analog of the vector product:
$$ (Y_1, Y_2, Y_3) \boxtimes (x_1,x_2,x_3) = (Y_2 x_3 -Y_3 x_2, 
Y_3 x_1 - Y_1 x_3, - Y_1 x_2 + Y_2 x_1)~.$$

In both cases, the fact that $b$ depends only on the equivalence class of $U'$
(under the action of trivial deformations) has an ``abstract''
interpretation. As defined above, for each polygon $p$ in $S^2$, $b$
defines a quadratic form on the space of isometric first-order deformations of
$p$, with values in $\R^3$. Let $l=(l_1, \cdots, l_n)$ be the edge lengths of
$p$. Since $b$
behaves ``well'' under the action of $SO(3)$ on $S^2$, $b$ actually defines a
quadratic form on $T_{[p]}\cP_S(l)$, where $[p]$ is the projection of $p$ in
$\cP_S(l)$, with values in the vector bundle
over $\cP_S(l)$ which is defined, from the trivial $\R^3$-bundle over
the space of polygons in $S^2$ with edge lengths given by $l$,
by taking the natural action of $SO(3)$ on both this space of polygons and
$\R^3$. 

\pg{A positivity property.}

The quadratic form $b$ defined above has a striking geometric
property when $p$ is a convex spherical or hyperbolic 
polygon. By a {\it convex polygon}, we mean a polygon which is the boundary of
a convex domain in $S^2$ (resp. $H^2$), which we call the {\it interior} of
$p$, and denote by $\mbox{int}(p)$.
We state this property
first in the spherical setting. Recall that, given a
convex spherical polygon $p$, the {\it dual polygon} $p^*$ is a convex polygon
whose interior is the set of points in $S^2$ which have positive scalar
product with the vertices of $p$ (see
e.g. \cite{coxeter-projective,coxeter-non-euclidean}). 
Each of its vertices is at distance $\pi/2$
of an edge of $p$, and conversely. 

\begin{thn}{B$_S$}
Let $p$ be a convex polygon in $S^2$, and let $U$ be 
a non-trivial infinitesimal
first-order deformation of $p$. Then $b(U)\in (\R_+\setminus\{ 0\})\mbox{int}
(p^*)$, i.e. $b(u)$ is contained in the positive cone over the interior of 
$p^*$.
\end{thn}

There is a geometric interpretation to this property. The space of polygons
with given edge lengths $l=(l_1, \cdots, l_n)$, $\cP_{S}(l)$, has a natural
map $\phi$ to $\R^n$, where $\phi(p)$ is the family of the angles of $p$. 
By Theorem
A$_S$, the image is a locally a submanifold except when $p$ has all its
vertices on a great circle. We will see in section 3 that $b$ is strongly
related to the second fundamental form of this submanifold, so that Theorem
B$_S$ translates as a convexity property: if $p$ is a convex polygon, the
image of $\phi$, in the neighborhood of $\phi(p)$, has a second fundamental
form which is positive definite in some directions. Theorem B$_S$ is related
to a result of Volkov \cite{volkov-polygones} on (non infinitesimal) isometric
deformations of spherical polygons. 

The same result holds in the hyperbolic or the de Sitter setting. Now the dual
of a hyperbolic polygon is a convex de Sitter polygon, and conversely; both
can be defined as in the sphere, using the Minkowski metric on $\R^3_1$ (see
e.g. \cite{coxeter-de-sitter,coxeter-non-euclidean}, where the ``polarities''
corresponding to the duality used here has a very central role).

\begin{thn}{B$_H$}
Let $p$ be a convex polygon in $H^2$
and let $U$ be a non-trivial infinitesimal
first-order deformation of $P$. Then $b(U)\in (\R_+\setminus \{ 0\})\mbox{int}
(p^*)$, i.e. $b(U)$ is in the positive cone over the interior of $p^*$.
\end{thn}

\pg{An isoperimetric statement for hyperbolic polygons.}

The first application of Theorems A$_H$ and B$_H$
that we consider is to an isoperimetric problem for hyperbolic
polygons. The next theorem is hyperbolic analog of a
statement which was proved by Steiner \cite{steiner} in the spherical setting,
but was known earlier in the Euclidean case (see \cite{siegel} for 
recent progress on more elaborate statements even in the Euclidean case). 

Given $l=(l_1, \cdots, l_n)$, we call $\cP_S^c(l)$ the space of convex
polygons in $S^2$ with edge lengths equal to the $l_i$, considered up to
global isometries in $S^2$, and $\cP_H^c(l)$ the
space of convex hyperbolic polygons with edge lengths equal to the $l_i$,
again considered up to global isometries.

\btm \label{tm:isoper-hyper}
Let $l=(l_1, \cdots, l_n)$ be such that $\cP^c_H(l)$ is not empty; then there
exists a unique element $p\in \cP^c_H(l)$ which has its vertices either on a
circle, on a horocycle, or on a connected set of points at fixed distance from
a geodesic. Moreover, $p$ has maximal area. 
\etm

There is a similar statement, which is slightly simpler, in the Euclidean
plane and in the sphere; in both cases, there is a unique polygon of maximal
area (among the  polygons with fixed edge lengths) and it has its vertices on
a circle. The spherical 
result was apparently discovered by Steiner \cite{steiner}, but he writes that
the Euclidean statement was previously known. The proof given here for
hyperbolic polygons also works in the spherical setting.

The first point of the proof is that, as a direct consequence of
Theorem A$_H$, a hyperbolic polygon is a
critical point of the area, restricted to polygons with the same edge lengths,
if and only if its vertices are either on a circle, a
horocycle, or a connected component of the space of points at fixed distance
from a geodesic. 

Theorem B$_S$ (resp. B$_H$) also has an interesting meaning
in this context: it implies that the area, as a function defined on the space
of convex polygon with given edge lengths, is ``often'' strictly concave for a
natural metric. 

The second point of the proof is that, on the boundary of the space of convex
polygons with given edge lengths, the interior normal derivative of the area
is positive. It follows that there exists at least one local maximum of the
area in the interior. 
But a direct argument shows that, given the edge lengths, there is at
most one polygon with its vertices on a circle, a horocycle, or on one
connected component of the set of points at fixed distance from a geodesic. 

\pg{The infinitesimal rigidity of convex polyhedra.}

The positivity property of $b$ in the spherical setting leads to a simple
proof of the infinitesimal rigidity of convex polyhedra in the Euclidean
space. Although the first proof was given by Dehn \cite{dehn-konvexer}, 
it follows from the ideas of Legendre \cite{legendre} and Cauchy
\cite{cauchy}. Other proofs have been given, in particular by Kann \cite{kann}
and Filliman \cite{filliman}.

\begin{thn}{C}{ \bf (Legendre, Cauchy, Dehn)}
Let $P$ be a convex polyhedron in $\R^3$. Any first-order deformation of $P$
which does not change its combinatorics or the metrics on its faces is
trivial, i.e. induced by a global Killing field. 
\end{thn}

The proof given here is in section 5. It bears some relations with a proof of
the global rigidity of convex Euclidean polyhedra discovered by Pogorelov
\cite{pogo-polygones}. 

\pg{Fuchsian polyhedral surfaces in the Minkowski space.}

The techniques described here have a natural application to another rigidity
problem, concerning polyhedral surfaces in the Minkowski 3-dimensional
space. We consider space-like polyhedra; rather than closed polyhedra,
which can be defined as images of convex polyhedral maps from the sphere, we
consider {\it equivariant polyhedra}, which are the images of a polyhedral map
from the universal cover of a surface of genus at least $2$ which is {\it
  equivariant} (more precise definitions can be found in section 6). We will
prove that, among those surfaces, those which are ``Fuchsian'' -- the
associated representation from the fundamental group of the surface to the
isometry group of $\R^3_1$ has its image in
$SO(2,1)$ -- are infinitesimally rigid, i.e. any first-order
deformation of those surfaces which does not change the induced metric is
trivial, this is Theorem \ref{tm:minkowski}.

There is an analogous result for smooth, equivariant, Fuchsian, convex
surfaces \cite{iie}, which was proved by related methods. Actually, finding a
polyhedral version of the results of \cite{iie} was the main motivation for
the present work, although the by-products turned out to be rather more
interesting. This question was also studied by I. Iskhakov \cite{iskhakov},
who provided some partial results. One feature of the proof of the rigidity of
convex polyhedra given here is the existence of a distinguished point, which
appears quite artificial in the Euclidian context; for surfaces in the
Minkowski space, however, it is quite natural and even necessary, since the
point which is fixed by the representation of the fundamental group already
plays a special role.

The infinitesimal rigidity of Fuchsian equivariant surfaces in the Minkowski
space is equivalent, thanks to the ``Pogorelov map'' used e.g. 
in \cite{iie}, to
similar rigidity statements in the de Sitter or anti-de Sitter space, so that
polyhedral rigidity results in those spaces could be obtained as a consequence
of the rigidity theorem for polyhedral surfaces in $\R^3_1$ proved
here. However, another proof of the rigidity of Fuchsian equivariant surfaces
in the de Sitter or the anti-de Sitter space has recently been developed
by F. Fillastre (in preparation), and it applies in particular in the
polyhedral setting. So it should also be possible to prove the result stated
here -- on polyhedral Fuchsian surfaces in $\R^3_1$ -- from the statements
obtained by Fillastre in the de Sitter or the anti-de Sitter space. 

To prove Theorem \ref{tm:minkowski}, we will follow the proof of Theorem C,
but the polygons that will be considered will be in the de Sitter plane rather
than in the sphere. 
Recall that, in addition to the hyperbolic plane, $\R^3_1$ contains another
quadric, the de Sitter plane $S^2_1$, which is a complete, constant curvature
2-dimensional Lorentz manifold (see
e.g. \cite{coxeter-de-sitter,coxeter-non-euclidean,O}): 
$$ S^2_1 = \{ (x,y,z)\in \R^3_1 ~ | ~ \langle (x,y,z),(x,y,z)\rangle =
1\}~. $$ 

Polygons in the de Sitter plane can be defined as in the sphere; the edges can
be of different types, either space-like, light-like or time-like. The notion
of angle is more subtle than in the sphere, see e.g. \cite{shu,cpt}, and it is
quite natural to consider the angles of a polygon as complex numbers,
with real part either $\pm \pi/2$ or $\pi$. A key point is that, with those
definitions and the corresponding definitions for the lengths of the edges,
the main triangle formulas in the sphere remain valid in the de Sitter plane,
a simple fact which we prove in section 6 for completeness. 

The first-order deformations of polygons in the de Sitter plane, in terms of
the first-order variation of the angles, can be described as in the sphere or
the hyperbolic plane. 

\begin{thn}{$\mbox{A}_{dS}$}
  Let $p=(v_1, \cdots, v_n)$ be a polygon in the de Sitter plane. 
  Let $\alpha_1, \cdots,
  \alpha_n$ be its angles, and let $\alphad_1, \cdots, \alphad_n\in \R$
  be a first-order variation of its angles induced by an isometric first-order
  deformation of $p$. Then:
$$ \sum_{i=1}^n \alphad_i v_i=0~, $$
where the $v_i$ are considered as points in $\R^3_1$.

Conversely, if this
equation is satisfied by an $n$-uple $(\alphad_1, \cdots,\alphad_n)$ and
moreover the $v_i$ are not all on geodesic, then there exists an
isometric first-order deformation of $p$ such that the $\alphad_i$ are the
associated first-order variations of the $\alpha_i$.
\end{thn}

In this context, the quadratic form $b$ can be defined as in the sphere. It
has the same positivity property for convex, space-like polygons (a notion
which is defined here with some care, see section 6). The notion of polygon
dual to a space-like, convex polygon is defined in section 6 as in the sphere,
it is a polygon in the hyperbolic plane. 

\begin{thn}{B$_{dS}$}
Let $p$ be a convex, space-like polygon in $S^2_1$
and let $U$ be an infinitesimal
first-order deformation of $P$. Then $b(U)\in (\R_+\setminus \{ 0\})\mbox{int}
(p^*)$, i.e. $b(U)$ is contained in the positive cone over the interior of the
dual polygon $p^*$.
\end{thn}

The argument used in the proof of Theorem C, translated to the Minkowski
setting, yields a rigidity statement for Fuchsian polyhedral surfaces in
$\R^3_1$, Theorem \ref{tm:minkowski}.

\pg{Metrics on moduli spaces.}

One consequence of Theorem $\mbox{B}_S$ is that it shows that the
invariant $b$ can be used to define some natural metrics on the moduli space
of convex polygons in $S^2$ (and similarly in the hyperbolic plane). This is
described in section 4, along with some relations to a ``natural'' metric on
Euclidean polygons with fixed angles, which has some interesting properties. 

\pg{Notations.} In all the paper, we set: $\R_+^*:=\R_+\setminus \{ 0\}$, and 
$\R_-^* := \R_-\setminus \{ 0\}$

\section{Deformations of polygons}

\pg{The moduli space of polygons.}

The geometry of the moduli space of polygons has been studied rather
extensively, in particular by Kapovich and Millson
(see e.g. \cite{kapovich-millson,kapovich-millson-sphere}). 
We recall here only some very
elementary properties which should clarify parts of the proofs below. We
consider here polygons in the plane, however all the comments in this
paragraph hold also for spherical, hyperbolic or de Sitter polygons, with some
obvious adaptations. 

Let $n\geq 3$, let $P_{n,E}$ be the space of polygons with $n$ vertices in
$\R^2$. Recall that a polygon is a sequence of vertices $v_1,\cdots, v_n=v_0$
such that, for each $i$, $v_i\neq v_{i+1}$. The isometry group $\Isom(\R^2)$
of $\R^2$ acts without fixed points on $P_{n,E}$, so we consider the quotient
$\cP_{n,E}=P_{n,E}/\Isom(\R^2)$. $\cP_{n,E}$ is a smooth manifold of dimension
$2n-3$, which can also be considered as a smooth algebraic variety. 

There is a family of $n$ functions naturally defined on $\cP_{n,E}$; if $p\in
P_{n,E}$ is a polygon, with $p=(v_1, \cdots, v_n)$, then:
$$ \lambda_i(p):=d(v_i,v_{i+1})~. $$
Those $n$ functions are clearly invariant under the action of $\Isom(\R^2)$, so
they define natural functions on $\cP_{n,E}$, which we still call
$\lambda_i$. Note that it would be algebraically more natural to consider the
{\it squares} of the distances, but we stick to the more natural definition
from an elementary geometry viewpoint (for spherical polygons one could
consider the cosine of the distance, and for hyperbolic polygons the cosh). 

Now let $(l_1, \cdots, l_n)$ be an $n$-uple of positive numbers, recall that
$\cP_{E}(l)$ is the (moduli) space of polygons with edge lengths $l_i$. So:
$$ \cP_E(l) := \{ p\in \cP_{n,E} ~ | ~ \forall i\in \{ 1,\cdots, n\},
\lambda_i(p) = l_i \}~, $$
so that $\cP_E(l)$ is an algebraic subvariety of $\cP_{n,E}$. 
Let $p\in \cP_E(l)$, consider the space of its first-order infinitesimal
deformations:
$$ \cT_p := \ker(d\lambda_1)\cap \cdots \cap \ker(d\lambda_n)\subset
T_p\cP_{n,E}~. $$ 
Clearly $\cT_p$ has dimension at least $n-3$, since it is defined by $n$
equations. If those $n$ equations are linearly independent, then, by the
inverse function theorem, in the
neighborhood of $p$, $\cP_E(l)$ is 
a smooth submanifold of $\cP_{n,E}$ (and also a
smooth algebraic subvariety) of dimension $n-3$. 

\pg{Euclidean polygons.}

We first indicate the proof of Proposition \ref{pr:euclidien}. 
The simplest proof is perhaps obtained 
by taking a limit, in Theorem A$_S$ below, when the lengths of the edges go
to $0$. It is however possible to give a simpler 
direct proof (this was pointed out by Sergiu Moroianu). 

Let $p=(v_1,
\cdots, v_n)$ be a Euclidean polygon, with edge lengths $l_1, \cdots, l_n$, so
that, for each $i\in \{ 1, \cdots, n\}$, $l_i=\| v_{i+1}-v_i\|\neq 0$ (with
$v_{n}=v_0$). For each $i$, let $\theta_i$ be the (oriented) angle between
the oriented $x$-axis and $v_{i+1}-v_i$. The possible values of the angles
$\theta_1, \cdots, \theta_n$, in the neighborhood of $p$,
are defined by condition that $ \sum_{i=1}^{n} l_i e^{i\theta_i} = 0 $,
because $v_{n+1}=v_1$. 

In an isometric first-order deformation of $p$, it follows that: 
$$ \sum_{i=1}^{n} \thetad_i l_i e^{i\theta_i} = 0~. $$
Suppose (adding a trivial deformation if necessary) 
that, in the deformation of $p$ we consider, the direction of
$v_1-v_0$ does not vary. Then, for each $i$, $\thetad_i =
\sum_{j=1}^i\alphad_j$. Therefore:
$$ 0 = \sum_{i=1}^{n} l_i e^{i\theta_i} \sum_{j=1}^i \alphad_i =
\sum_{j=1}^{n} \alphad_j\sum_{i=j}^{n} l_i e^{i\theta_i}~. $$
We can also suppose -- still 
without loss of generality -- that $v_0=0$; then, for each $j$: 
$$ v_j=\sum_{i=0}^{j-1} l_ie^{i\theta_i} = - \sum_{i=j}^{n} l_i
e^{i\theta_i}~, $$
and it follows that:
$$ \sum_{j=1}^{n} \alphad_j v_j=0~. $$

In addition, it is well known that, for any Euclidean polygon:
$$ \sum_{i=1}^n \alpha_i = 2\pi k~, $$
where $k$ is the winding number of the polygon. It follows that, in any
first-order deformation, $\sum_{i=1}^n \alphad_i=0$.

The two conditions in the statement of Proposition \ref{pr:euclidien} are
linearly independent if and only if the $v_i$ are not all on a line. 
When the $v_i$ are not collinear, 
the vector space of $n$-uples $(\alphad_1, \cdots,
\alphad_i)$ satisfying them has dimension $n-3$. Therefore, the arguments in
the previous paragraph show that 
this space has dimension exactly $n-3$, and the possible first-order
variations of the $\alpha_i$ are exactly the solutions of the two
equations in the statement of the proposition.

\pg{Spherical polygons.} 

We now prove Theorem A$_S$. A direct proof is possible\footnote{One such proof
  is available in the first version of this paper, see
  \texttt{http://front.math.ucdavis.edu/math.DG/0410058} and then choose
  ``v1''.}, but the proof given here is geometric and much simpler. 

Let $p=(v_1, \cdots, v_n)$ be a spherical
polygon, and let $U$ be a first-order isometric deformation of $p$. Consider
the polyhedral cone $C$ over $p$, i.e. the union of the half-lines with
endpoint at $0$ which intersect $p\subset S^2\subset \R^3$. $C$ has $n$ faces
$f_1, \cdots, f_n$, with $f_i$ equal to the cone over the edge $e_i$ of
$p$. Since $\alpha_i$ is the angle between the edges $e_i$ and $e_{i+1}$ of
$p$, it is also equal to the dihedral angle between $f_i$ and $f_{i+1}$.

The first-order deformation $U$ of $P$ determines a first-order deformation
$V$ of $C$. For each $i\in \{ 1,\cdots, n\}$, $V$ acts on each face of $C$ as
a Killing field, i.e. there exist vectors $u_i$ and $y_i$ such that $V$ acts
on $f_i$ as the vector field $\kappa_i$ defined by:
$$ \forall x\in \R^3, \kappa_i(x) = u_i + y_i\times x~. $$
Then $\kappa_{i+1}-\kappa_i=0$ on the intersection of $f_i$ and $f_{i+1}$, and
the difference between $\kappa_i$ and $\kappa_{i+1}$ is equal to an
infinitesimal rotation of angle $d\alpha_i(U)$ and axis $f_i\cap f_{i+1}$. It
follows that: 
$$ y_{i+1} = y_i + d\alpha_i(U) v_i~, $$
and the result follows since the sum of the differences $y_{i+1} - y_i$, for
$i$ going from $1$ to $n$, vanishes.

\pg{Hyperbolic polygons.}

The same argument can be used to prove Theorem A$_H$. $C$ is now a polyhedral
cone in the Minkowski space $\R^3_1$, and its faces are time-like. The
restriction of the infinitesimal deformation $V$ to each face $f_i$ is a
Killing vector, which is now of the form: $\kappa_i(x)=u_i+y_i\boxtimes
x$. The proof then proceeds as for the sphere.

\section{A positivity result}

In this section we prove Theorem B$_S$, as well as its hyperbolic analog
Theorem B$_H$, after some preliminary computations concerning the first-order
deformations of quadrilaterals.

\pg{Deformations of quadrilaterals.}

\begin{figure}[h] \label{fg:quadri}
\centerline{\psfig{figure=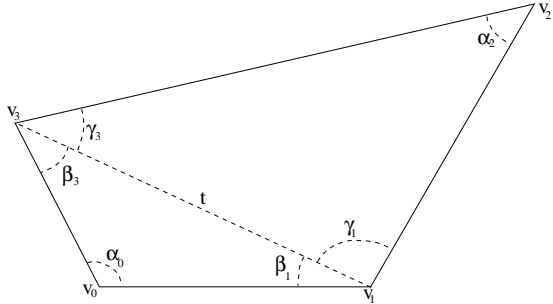,height=4cm}}
\caption{Deformation of a non-degenerate quadrilateral.}
\end{figure}

We consider here a quadrilateral $q=(v_1, v_2, v_3, v_4)$ 
for which no 3 of the
vertices are collinear, and set $v_0:=v_4$. 
We use the notations apparent in Figure
\ref{fg:quadri}, in particular $t$ is the distance between $v_1$ and $v_3$. We
also call $d_{i,j}$ the distance between $v_i$ and $v_j$. 

Consider a first-order isometric deformation of $q$, suppose for instance that
$v_0$ and $v_1$ are fixed (this can be achieved by adding a trivial
deformation). Then $v_2$ moves in the direction of the circle of center $v_1$
containing it. Since $v_1, v_2$ and $v_3$ are not collinear, the circle of
center $v_3$ containing $v_2$ is not tangent to the circle of center $v_1$
containing $v_2$, so that any first-order displacement of $v_2$ induces a
non-zero first-order variation of the distance from $v_2$ to $v_3$.

The same argument shows that any first-order displacement of $v_3$, preserving
its distance to $v_0$, induces a non-zero first-order variation of its
distance to $v_2$. It follows that there is a 1-dimensional vector space of
isometric 
first-order deformations of $q$ (up to the trivial deformations), they are
parametrized for instance by $t$, the distance from $v_1$ to $v_3$. So we
consider the first-order deformation $U$ of $q$ such that $dt(U)=1$.

To simplify notations, we call $d_{ij}$ the distance between $v_i$ and $v_j$,
for $0\leq i,j\leq 4$. Then a well-known spherical geometry formula states
that: 
\beq \label{eq:cos}
\cos(t) = \cos(d_{01}) \cos(d_{03}) + \cos(\alpha_0) \sin(d_{01})
\sin(d_{03})~, \eeq
so that:
\beq \label{eq:alpha0}
d\alpha_0(U) = \frac{\sin(t)}{\sin(d_{01}) \sin(d_{03})
  \sin(\alpha_0)}~. \eeq
The same computation (or a symmetry argument) shows that:
\beq \label{eq:alpha2}
d\alpha_2(U) = \frac{\sin(t)}{\sin(d_{12}) \sin(d_{23})
  \sin(\alpha_2)}~. \eeq
Moreover, equation (\ref{eq:cos}), applied to the triangle $(v_0,v_1, v_3)$, 
yields:
$$ \cos(\beta_1) = \frac{\cos(d_{03}) - \cos(d_{01})\cos(t)}{\sin(d_{01})
  \sin(t)}~, $$
so that:
\begin{eqnarray*}
- \sin(\beta_1) d\beta_1(U) & = &
\frac{\cos(d_{01})\sin^2(t) - (\cos(d_{03}) -
  \cos(d_{01})\cos(t))\cos(t)}{\sin(d_{01})\sin^2(t)} \\ 
& = &\frac{\cos(d_{01}) - \cos(d_{03})\cos(t)}{\sin(d_{01}) \sin^2(t)}~.
\end{eqnarray*}
But, by the sine formula for spherical triangles:
$$  \frac{\sin(\beta_1)}{\sin(d_{03})} = \frac{\sin(\alpha_0)}{\sin(t)}~, $$
so that: 
$$ d\beta_1(U) =
\frac{\cos(d_{03})\cos(t)-\cos(d_{01})}{\sin(d_{01})\sin(d_{03})
  \sin(t)\sin(\alpha_0)}~. $$
The same computation (or a symmetry argument) shows that:
$$ d\gamma_1(U) =
\frac{\cos(d_{23})\cos(t)-\cos(d_{12})}{\sin(d_{12})\sin(d_{23})
  \sin(t)\sin(\alpha_2)}~, $$
and, taking the sum, we obtain that:
\beq \label{eq:alpha1} d\alpha_1(U) =
  \frac{\cos(d_{03})\cos(t)-\cos(d_{01})}{\sin(d_{01})\sin(d_{03}) 
  \sin(t)\sin(\alpha_0)} + 
  \frac{\cos(d_{23})\cos(t)-\cos(d_{12})}{\sin(d_{12})\sin(d_{23})
  \sin(t)\sin(\alpha_2)}~. \eeq
By the same computation (or a symmetry argument):
\beq \label{eq:alpha3} 
d\alpha_3(U) =
  \frac{\cos(d_{01})\cos(t)-\cos(d_{03})}{\sin(d_{01})\sin(d_{03}) 
  \sin(t)\sin(\alpha_0)} + 
  \frac{\cos(d_{12})\cos(t)-\cos(d_{23})}{\sin(d_{12})\sin(d_{23})
  \sin(t)\sin(\alpha_2)}~. \eeq

\pg{Proof of Theorem B$_S$.} Let $p=(v_1, \cdots, v_n)$ be a convex spherical
polygon. Recall that, by Definition \ref{df:bUS}, for any
first-order deformation $U$ of $p$, defined up to the Killing fields, we have:
$$ b(U) = \sum_{i=0}^n d\alpha_i(U') dv_i(U')~, $$
where $U'$ is any representative of $U$, i.e. any first-order deformation of
$p$ corresponding to $U$ under the quotient by the trivial deformations. 

So $b$ is a quadratic form on $T_p\cP_S(l)$, where $l=(l_1, \cdots, l_n)$
is the $n$-uple of the edge lengths of $p$. It is natural to define a bilinear
form associated to $b$, which we call $b_2$; it is defined as follows: if
$U'_1$ and $U'_2$ are two first-order deformations of $p$, then:
$$ b_2(U'_1, U'_2) := \frac{1}{2}\left(\sum_{i=1}^n d\alpha_i(U'_1) dv_i(U'_2)
+ \sum_{i=1}^n d\alpha_i(U'_2) dv_i(U'_1)\right)~. $$
As for $b$, an important point is that if one adds a Killing field to $U'_1$
or $U'_2$, the result does not change. Moreover, it will be useful below to
note that each of the two sums in the definition of $b_2$ is invariant under
this transformation. Indeed, if $Y_1, Y_2\in \R^3$ are two vectors, let $V_1,
V_2$ be the trivial deformations defined by $dv_i(V_1)=Y_1\times v_i,
dv_i(V_2) = Y_2\times v_i$; then, for
each $i\in \{ 1,\cdots, n\}$, $d\alpha_i(U'_1+V_1)=d\alpha_i(U'_i)$, and:
\begin{eqnarray*}
\sum_{i=1}^n d\alpha_i(U'_1+V_1) dv_i(U'_2+V_2) & = &
\sum_{i=1}^n d\alpha_i(U'_1) (dv_i(U'_2) + Y_2\times v_i) \\
& = & \sum_{i=1}^n d\alpha_i(U'_1)dv_i(U'_2) + Y_2 \times \sum_{i=1}^n
d\alpha_i(U'_1)v_i \\ 
& = &
\sum_{i=1}^n d\alpha_i(U'_1)dv_i(U'_2) ~, 
\end{eqnarray*}
and the same computation can be applied to the second sum in the definition of
$b_2$. 

\blm \label{lm:decomposition}
Let $p=(v_1, \cdots, v_n=v_0)$ be a polygon such that no 3 vertices are
collinear. Let $U$ be an isometric first-order deformation of $p$, such that
$v_1$ and $v_2$ are fixed. There exists a unique decomposition:
$$ U = U_2 + U_3 + \cdots + U_{n-2} $$
such that, for all $i\in \{ 2, \cdots, n-2\}$, $U_i$ vanishes on $v_1, v_2,
\cdots, v_i$, and acts on $v_{i+1}, \cdots, v_{n-1}, v_n=v_0$ as a
rigid motion. 
\elm

\begin{proof}
We define a sequence of isometric first-order deformations of $p$ recursively,
as follows. First note that, since $v_0, v_1$ and $v_3$ are not collinear, 
there is a unique first-order isometric
deformation $V_2$ of the quadrilateral $(v_0, v_1, v_2, v_3)$ such that $v_1$
and $v_2$ are fixed, and that the first-order displacement of $v_3$ is the
same for $V_2$ and for $U$. Then define a first-order deformation $U_2$ as
follows: 
\begin{itemize}
\item $v_1$ and $v_2$ are fixed.
\item $v_3$ moves under $U_2$ as under $V_2$, i.e. as under $U$.
\item The restriction of $U_2$ to $v_3, v_4, \cdots, v_n=v_0$ is equal to the
  restriction to those vertices of a Killing field, i.e. the corresponding
  part of $p$ moves in a ``rigid'' way. 
\end{itemize}
It is possible to define such a deformation because, at first order, the
distance between $v_3$ and $v_0$ does not vary under the deformation $V_2$.

Now consider the first-order deformation $U-U_2$ of $p$. It is isometric ---
as the difference of two isometric deformations --- and it vanishes at
$v_1, v_2$ and $v_3$. This shows, because $v_1, v_4$ and $v_0$ are not
collinear, that there exists a unique
isometric first-order deformation $V_3$ of the quadrilateral 
$(v_1, v_3, v_4, v_0)$ which vanishes at $v_1$ and $v_3$ and acts on $v_4$ as
$U-U_2$. Let $U_3$ be the unique first-order deformation of $p$ which vanishes
at $v_1, v_2$ and $v_3$, and which acts on $v_4, \cdots, v_{n-1}, v_0$ as
a Killing field  (i.e. ``rigidly''). 

Now $U-U_2-U_3$ is an isometric
first-order deformation of $p$ which vanishes at $v_1, v_2, v_3$ and $v_4$, and
we can iterate this construction until we obtain a decomposition:
$$ U = U_2 + \cdots + U_{n-2}~, $$
which has the property described in the lemma.
\end{proof}

Consider now a first-order deformation $U$ of $p$. Since $p$ is convex, it is
non-degenerate in the sense of Lemma \ref{lm:decomposition}. Applying this
lemma yields a decomposition $U = U_2+ \cdots +U_{n-2}$ such that, 
for all $i\in \{ 2, \cdots, n-2\}$, $U_i$ vanishes on $v_1, v_2,
\cdots, v_i$, and acts on $v_{i+1}, \cdots, v_{n-1}, v_n=v_0$ as a
rigid motion. 

\blm \label{lm:diagonal}
For each $i,j\in \{ 2, \cdots, n-2\}$ with $i\neq j$, $\langle b_2(U_i,
U_j), v_1\rangle=0$. Therefore: 
$$ \langle b(U), v_1\rangle = \left\langle \sum_{i=1}^n b(U_i),
  v_1\right\rangle ~. $$
\elm

\bpv
Suppose (without loss of generality) that $i<j$. 
For each $i\in \{ 2,\cdots, n-2\}$, we add a Killing field to
$U_i$ to obtain a first-order deformation $U'_i$ which acts on $v_1, v_2,
\cdots, v_i$ as a rigid motion, and vanishes on $v_{i+1}, \cdots, v_{n-1},
v_n=v_0$. It follows from the remark above that:
$$ b_2(U_i, U_j) = \frac{1}{2}\left(\sum_{k=1}^n d\alpha_k(U_i) dv_k(U_j)
+ \sum_{i=1}^n d\alpha_k(U'_j) dv_k(U'_i)\right)~. $$
Note that $d\alpha_k(U_i)$
vanishes for all values of $k$ except $0, 1, i$ and $i+1$. But $dv_k(U_j)$
vanishes, by definition of $U_j$, for $k=1, \cdots, j$. Symmetrically,
$d\alpha_k(U'_j)$ vanishes except for $k=0,1,j,j+1$, while $dv_k(U'_i)$
vanishes for $k=i+1, \cdots, n$, so that:
$$ b_2(U_i, U_j) = \frac{1}{2} (d\alpha_0(U_i)dv_0(U_j) + d\alpha_1(U'_j)
dv_1(U'_i))~. $$
Since $U_j$ is an isometric first-order deformation, we have:
$$ 0 = d\langle v_0,v_1\rangle (U_j) = \langle dv_0(U_j), v_1\rangle + \langle
v_0, dv_1(U_j)\rangle~. $$
But $dv_1(U_j)=0$ by definition of $U_j$, so that $\langle dv_0(U_j),
v_1\rangle =0$. Since $dv_1(U'_i)\in T_{v_1}S^2$,  
$\langle v_1,dv_1(U'_i)\rangle
=0$, and it follows that $\langle b_2(U_i,U_j),v_1\rangle =0$.  
\epv

The value of $\langle b(U_i), v_1\rangle$ is given by the next lemma.

\blm \label{lm:calcul}
Let $q:=(v_1, v_i, v_{i+1}, v_0)$, and let $\alpha'_1, \alpha'_i,
\alpha'_{i+1}, \alpha'_0$ be its angles. Let $V$ be the first-order isometric
deformation of $q$ which vanishes on 
$v_1$ and $v_i$ and under which the distance between $v_1$ and $v_{i+1}$
varies at speed $1$. Then, if $d_{j,k}:=d(v_j, v_k)$, we have:
\beq \label{eq:calcul}
\langle v_1, d\alpha'_{i+1}(V) dv_{i+1}(V) + d\alpha'_0(V)dv_0(V)\rangle =
\frac{\sin^2(d_{1,i+1})\sin(\alpha'_1)}{\sin(d_{0,i+1}) \sin(d_{i,i+1})
  \sin(\alpha'_0)\sin(\alpha'_i)}~. \eeq
\elm

\begin{proof}
We have already noted that $\langle v_1,dv_0(V)\rangle=0$ because $V$ does not
change, at first order, the distance between $v_0$ and $v_1$. So the only
non-vanishing term is the one involving $\langle v_1, dv_{i+1}(V)\rangle$,
and: 
$$ \langle v_1, dv_{i+1}(V)\rangle = d\langle v_1, v_{i+1}\rangle(V) =
d\cos(d_{1,i+1})(V) = - \sin(d_{1,i+1})~. $$
But $d\alpha_{i+1}(V)$ is given by equation (\ref{eq:alpha3}); we now have
slightly different notations, $U$ is replaced by $V$, the indices $2$
(resp. $3$) by $i$ (resp. $i+1$). So:
$$ d\alpha'_{i+1}(V) =
\frac{\cos(d_{0,1})\cos(t)-\cos(d_{0,i+1})}{\sin(d_{0,1})\sin(d_{0,i+1})  
  \sin(t)\sin(\alpha'_0)} + 
  \frac{\cos(d_{1,i})\cos(t)-\cos(d_{i,i+1})}{\sin(d_{1,i})\sin(d_{i,i+1})
  \sin(t)\sin(\alpha'_i)}~, $$
where $t:=d_{1,i+1}$. Following the notations of Figure \ref{fg:quadri}, 
we call $\beta'_1$
and $\beta'_{i+1}$ the angles at $v_1$ and $v_{i+1}$, respectively, of the
triangle $(v_0, v_1, v_{i+1})$, and we call $\gamma'_1$ and $\gamma'_{i+1}$
the angles at $v_1$ and $v_{i+1}$, respectively, 
of the triangle $(v_1, v_i, v_{i+1})$. Using equation (\ref{eq:cos}), we get:
\begin{eqnarray*}
d\alpha'_{i+1}(V) & = &
- \frac{\cos(\beta'_1)\sin(d_{0,1})\sin(t)}{\sin(d_{0,1})\sin(d_{0,i+1})  
  \sin(t)\sin(\alpha'_0)} -
  \frac{\cos(\gamma'_1)\sin(d_{1,i})\sin(t)}{\sin(d_{1,i})\sin(d_{i,i+1})
  \sin(t)\sin(\alpha'_i)} \\ 
& = & - \sin(t) \left(\frac{\cos(\beta'_1)}{\sin(d_{0,i+1})
    \sin(\alpha'_0)\sin(t)} +
  \frac{\cos(\gamma'_1)}{\sin(d_{i,i+1})\sin(\alpha'_i)\sin(t)} \right) \\ 
& = & - \frac{\sin(t)}{\sin(d_{0,i+1}) \sin(d_{i,i+1})
    \sin(\alpha'_0) \sin(\alpha'_i)} \times \\
& \times & \left( \cos(\beta'_1)
    \frac{\sin(d_{i,i+1})\sin(\alpha'_i)}{\sin(t)} + \cos(\gamma'_1)
  \frac{\sin(d_{0,i+1})
    \sin(\alpha'_0)}{\sin(t)}\right)~.
\end{eqnarray*}
But, by the sine law for spherical triangles:
$$ \frac{\sin(\alpha'_i)}{\sin(t)} = \frac{\sin(\gamma'_1)}{\sin(d_{i,i+1})}~,
~~ \frac{\sin(\alpha'_0)}{\sin(t)}=\frac{\sin(\beta'_1)}{\sin(d_{0,i+1})}~, $$ 
so that:
\begin{eqnarray*}
d\alpha'_{i+1}(V) & = & - \sin(t) \frac{\cos(\beta'_1) \sin(\gamma'_1) +
    \cos(\gamma'_1) \sin(\beta'_1)}{\sin(d_{0,i+1}) \sin(d_{i,i+1}) 
    \sin(\alpha'_0) \sin(\alpha'_i)} \\
& = & - \sin(t) \frac{\sin(\alpha'_1)}{\sin(d_{0,i+1}) \sin(d_{i,i+1}) 
    \sin(\alpha'_0) \sin(\alpha'_i)}~,
\end{eqnarray*}
and the result follows.
\end{proof}

It follows directly from Lemma \ref{lm:diagonal} and from Lemma
\ref{lm:calcul} that $\langle b(U), v_1\rangle$ is a sum of positive terms, so
it is positive. By symmetry the same holds for all the other vertices. Since
$b(U)$ has positive scalar product with all the vertices of $p$, it is
contained in the positive cone over
the interior of the dual polygon $p^*$, and this proves Theorem B$_S$.

\pg{Proof of Theorem B$_H$.}

The same proof applies in the hyperbolic case, with some small differences in
the computations. We only state the hyperbolic analog of Lemma
\ref{lm:calcul}.

\blm \label{lm:calcul-h}
Let $q:=(v_1, v_i, v_{i+1}, v_0)$, 
and let $\alpha'_1, \alpha'_i, \alpha'_{i+1},
\alpha'_0$ be its angles. Let $V$ be the first-order isometric
deformation of $q$ which vanishes on 
$v_1$ and $v_i$ and under which the distance between $v_1$ and $v_{i+1}$
varies at speed $1$. Then, if $d_{j,k}:=d(v_j, v_k)$, we have:
$$ \langle v_1, d\alpha'_{i+1}(V) dv_{i+1}(V) + d\alpha'_0(V)dv_0(V)\rangle =
\frac{\sinh^2(d_{1,i+1})\sin(\alpha'_1)}{\sinh(d_{0,i+1}) \sinh(d_{i,i+1})
  \sin(\alpha'_0)\sin(\alpha'_i)}~. $$
\elm

The proof is based on some computations which are quite parallel to those made
above for first-order deformations of spherical quadrilaterals, but now for
hyperbolic quadrilaterals.
We use again the notations apparent in Figure \ref{fg:quadri}.  
The basic triangle formula is now:
\beq \label{eq:cosh} 
\cosh(t) = \cosh(d_{01})\cosh(d_{03}) - \cos(\alpha_0)
\sinh(d_{01})\sinh(d_{03})~, \eeq
from which it follows that, in first-order deformation $U$ such that $t$
varies at speed $1$, we have:
$$ d\alpha_0(U) = \frac{\sinh(t)}{\sin(\alpha_0) \sinh(d_{01})
  \sinh(d_{03})}~, $$
and the same computation shows that: 
$$ d\alpha_2(U) = \frac{\sinh(t)}{\sin(\alpha_2) \sinh(d_{12})
  \sinh(d_{23})}~. $$
Moreover, equation (\ref{eq:cosh}) yields:
$$ \cos(\beta_1) = \frac{- \cosh(d_{03}) + \cosh(d_{01})\cosh(t)}{\sinh(d_{01})
  \sinh(t)}~, $$
so that:
\begin{eqnarray*}
- \sin(\beta_1) d\beta_1(U) & = &
\frac{\cosh(d_{01})\sinh^2(t) - ( - \cosh(d_{03}) +
  \cosh(d_{01})\cosh(t))\cosh(t)}{\sinh(d_{01})\sinh^2(t)} \\ 
& = &\frac{- \cosh(d_{01}) + \cosh(d_{03})\cosh(t)}{\sinh(d_{01}) \sinh^2(t)}~.
\end{eqnarray*}
But:
$$  \frac{\sin(\beta_1)}{\sinh(d_{03})} = \frac{\sin(\alpha_0)}{\sinh(t)}~, $$
so that: 
$$ d\beta_1(U) =
\frac{- \cosh(d_{03})\cosh(t) + \cosh(d_{01})}{\sinh(d_{01})\sinh(d_{03})
  \sinh(t)\sin(\alpha_0)}~. $$
By symmetry:
$$ d\gamma_1(U) =
\frac{- \cosh(d_{23})\cosh(t) + \cosh(d_{12})}{\sinh(d_{12})\sinh(d_{23})
  \sinh(t)\sin(\alpha_2)}~, $$
and, taking, the sum:
$$ d\alpha_1(U) =
  \frac{- \cosh(d_{03})\cosh(t) + \cosh(d_{01})}{\sinh(d_{01})\sinh(d_{03}) 
  \sinh(t)\sin(\alpha_0)} + 
  \frac{- \cosh(d_{23})\cosh(t) + \cosh(d_{12})}{\sinh(d_{12})\sinh(d_{23})
  \sinh(t)\sin(\alpha_2)}~. $$
By the same computation (or a symmetry argument):
$$ d\alpha_3(U) =
  \frac{ - \cosh(d_{01})\cosh(t) + \cosh(d_{03})}{\sinh(d_{01})\sinh(d_{03}) 
  \sinh(t)\sin(\alpha_0)} + 
  \frac{ - \cosh(d_{12})\cosh(t) + \cosh(d_{23})}{\sinh(d_{12})\sinh(d_{23})
  \sinh(t)\sin(\alpha_2)}~. $$

The end of the proof then proceeds as in the spherical case, we leave the
details to the interested reader.

\pg{The submanifold of angles for given edge lengths.}

One possible way to interpret geometrically the positivity property of $b$ for
spherical polygons is
in terms of a kind of convexity property of a submanifold of codimension $3$
in $\R^n$, defined as the set of possible angles of convex polygons with given
edge lengths. This submanifold is often smooth, and when it has singularities,
they are located at precise points. Its tangent space at each point is given
by Theorem A$_S$, while Theorem B$_S$ indicates that its second fundamental
form is always on ``one side'', as explained below. 

We consider here a spherical polygon $p$ with $n$ vertices, and call
$l=(l_1, \cdots, l_n)$ the family of its edge length. Then we call $\cA_S(l)$
the space of $n$-uples $(\alpha_1, \cdots, \alpha_n)$ such that there exists a
polygon $p'\in \cP_S(l)$ with angles equal to the $\alpha_i$. This defines an
map $\Phi_l$ from $\cP_S(l)$ into $\R^n$, which is with image 
$\cA_S(l)$, because a
spherical polygon is entirely defined -- up to global isometries in $S^2$ --
by its edge lengths and its angles.

Theorem A$_S$,
along with the remarks at the beginning of section 2, show that $\cA_S(l)$ is
locally a smooth submanifold of codimension $3$ of $(\R/2\pi\Z)^n$, 
except when all the
vertices of the polygon lie on a spherical geodesic. Clearly this is possible
only when there exists $(\epsilon_1, \cdots, \epsilon_n)\in \{ -1,1\}$ and
$k\in \Z$ such that:
$$ \sum_{i=1}^n \epsilon_i l_i=2k\pi~. $$
Then the singular points of $\cA_S(l)$ can happen only at the points of
$(\R/2\pi\Z)^n$ which have all their coordinates equal to $0$ or $\pi$.

Let $q\in P_S(l)$, with vertices $v_1, \cdots, v_n$. Let $\alpha_1, \cdots,
\alpha_n$ be the angles of $q$, we suppose that they are not all equal to $0$
or $\pi$, so that $\cA_S(l)$ is locally smooth in the neighborhood of
$\alpha:=(\alpha_1,\cdots, \alpha_n)$. 
The tangent space of $\cA_S(l)$ is described
by Theorem A$_S$. More precisely, there is a natural isomorphism between
$\R^3$ and the normal space of $\cA_S(l)$ at $\alpha$, defined as follows:
$$ \forall w\in \R^3, \Psi_{q}(w) := (\langle v_1, w\rangle, \cdots, \langle
v_n,w\rangle)~. $$
Theorem A$_S$ shows that $\Psi_{q}(\R^3)$ is the orthogonal to the tangent
space of $\cA_S(l)$ at $\alpha$, i.e. the normal space of $\cA_S(l)$ at
$\alpha$.

Now $b$ appears naturally as the second fundamental form of $\cA_S(l)$.

\brk \label{rk:II}
Let $U\in T_qP_S(l)$ be a first-order isometric deformation of $q$, and let
$w\in \R^3$. then:
$$ \langle \II(d\Phi_l(U), d\Phi_l(U)), \Psi_q(w)\rangle = - \langle
b(U),w\rangle ~, $$
where $\II$ is the second fundamental form of $\cA_S(l)$.
\erk

\bpv
Let $w\in \R^3$. Consider any extension of $U$ as a vector field tangent to
$P_S(l)$ in the neighborhood of $q$, and call $\nabla$ the flat connection on
$\R^n$. 
By definition of the second fundamental form of $\cA_S(l)$, we have:
\begin{eqnarray*}
\langle \II(d\Phi_l(U),d\Phi_l(U)), \Psi_q(w)\rangle 
& = & - \langle \nabla_{d\Phi_l(U)} \Psi_q(w), d\Phi_l(U)\rangle \\
& = & - \langle (\langle dv_i(U), w\rangle )_{i=1}^n, 
(d\alpha_i(U))_{i=1}^n\rangle \\
& = & - \sum_{i=1}^n \langle d\alpha_i(U) dv_i(U) ,w\rangle \\
& = & - \langle b(U), w\rangle~.  
\end{eqnarray*}
\epv

There is a direct translation of Theorem B$_S$ in this context, indicating a
kind of convexity property of $\cA_S(l)$ at the points which are the images of
convex polygons. 

\begin{prop}
Suppose that $q$ is convex. Let $Q^*$ be the polyhedral cone in
$N_{\phi_l(q)}\cA_S(l)$ 
which is dual to the cone $\Psi_q((\R_+^*) q)$, for the metric
induced on $N_q\cA_S(l)$ by the metric on $\R^n$. Then:
$$ \forall V\in T_{\Psi_l(q)}\cA_S(l), \II(V,V)\in - Q^*~. $$
\end{prop}

\bpv
By the previous remark and Theorem B$_S$:
$$ \forall U\in T_qP_S(l), \forall w\in \R_+^* q, \langle \II(d\Phi_l(U),
d\Phi_l(U)), \Psi_q(w)\rangle < 0~, $$
so that:
$$ \forall V\in T_{\Phi_l(q)}\cA_S(l), \forall W\in \Psi_q(\R_+^* q), \langle
\II(V,V), W\rangle < 0~, $$
and the result follows. 
\epv

The same considerations hold also for hyperbolic polygons, based on Theorem
A$_H$ for the description of the tangent space, and on Theorem B$_H$ for the
convexity property of the submanifold of the angles of polygons with given
edge lengths. 

In the Euclidean case, one can again consider the space $\cA_E(l)$ of angles
of polygons with edge lengths given by $l$. Then $\cA_E(l)$ is a codimension
$2$ submanifold of the hyperplane of equation: $\sum_{i=1}^n \alpha_i = 0$. 
However it does not appear to have any obvious ``convexity'' property.

\section{Isoperimetric problems for spherical or hyperbolic polygons}

The theorems given above -- Theorem A$_S$ and Theorem B$_S$ for spherical
polygons, Theorem A$_H$ and Theoorem B$_H$ for hyperbolic polygons -- have
very simple applications to isoperimetric questions in the two settings. The
spherical result was apparently discovered by Steiner \cite{steiner}, and is
analoguous to an Euclidean
statement which was known earlier (see e.g. \cite{siegel} for
recent progress on related but more elaborate Euclidean statements). 
The proofs are strongly related
to Theorems A$_S$ and B$_S$ (resp. A$_H$ and B$_H$).

\pg{Critical points of the area functional.}

First, Theorem A$_S$ and Theorem A$_H$ lead directly to a description of the
critical points of the area function over the space of polygons with given
edge lengths.
 
\blm \label{lm:critic-sphere}
Let $p$ be a spherical polygon, with edge lengths given by
$l=(l_1,\cdots, l_n)$. Suppose that the vertices
of $p$ are not all on a geodesic. Then $p$ is a critical point of the area in
$\cP_S(l)$  if and only if all vertices of $p$ are on a
circle.
\elm

\blm \label{lm:critic-hyper}
Let $p$ be a hyperbolic polygon, with edge lengths given by
$l=(l_1,\cdots, l_n)$. Suppose that the vertices
of $p$ are not all on a geodesic. Then $p$ is a critical point of the area in
$\cP_H(l)$  if and only if all vertices of $p$ are either on a circle, on a
horocycle, or on a connected component of the set of points at fixed distance
from a geodesic. 
\elm

Note that the polygons are not required to be convex. It follows from the
results given below that, for convex polygons, the critical points are
actually maxima, but it is probably not the case for non-convex polygons. 
Both lemmas are proved together. 

\bpv
By the Gauss-Bonnet formula, $p$ is a critical point of the area functional if
and only if, for any first-order isometric deformation of $p$, the sum of the
angles of $p$ remains constant. By Theorem A$_S$ (resp. Theorem A$_H$) this
happens if and only if:
$$ \forall \alphad_1,\cdots, \alphad_n\in \R, \left(
\sum_{i=1}^n \alphad_i v_i =0
\Rightarrow \sum_{i=1}^n \alphad_i=0\right)~. $$
By an elementary linear algebra argument using the transposed equation, 
this is true if and only if there
exists a vector $u\in \R^3$ (resp. $\R^3_1$) such that:
$$ \forall i\in \{ 1, \cdots, n\}, \langle u,v_i\rangle =1~, $$
which holds if and only if the $v_i$ are on an affine plane in $\R^3$
(resp. $\R^3_1$) not containing the origin. 
The result follows, because the intersection of $S^2$ with an affine
plane in $\R^3$ is a circle (or a point) while the intersection of $H^2$ with
a plane in $\R^3_1$ is either a circle, a horocycle, or a connected component
of the set of points at fixed distance from a geodesic.
\epv

\pg{Uniqueness of critical polygons.}

In light of the previous paragraph, it is interesting to remark that, for a
given set of edge lengths, there is at most one convex polygon which is a
critical point of the area. We first consider spherical polygons, for which
the result has been well-known for many years.

\blm \label{lm:unique-sphere}
Let $l_1, \cdots, l_n\in \R_+^*$. There is at most one convex polygon $p\in
\cP^c_S(l)$ which has its vertices on a circle. 
\elm

\bpv
We suppose that $l_1$ is the largest of the $l_i$. 

Suppose that the vertices of $p$ are on a circle of radius $l$ and center
$x_0$ in $S^2$. Then the edge of $p$ of length $l_i$ is ``seen'' from $x_0$
under an angle $\theta_i$, and the sine law for spherical triangles shows
that: 
$$ \sin(\theta_i/2) = \frac{\sin(l_i/2)}{\sin(l)}~. $$
Setting $S_i:=\sin(l_i/2)$ and $S:=\sin(l)$, we have again two
possibilities: 
\begin{itemize}
\item Either $x_0$ is in the interior of $p$, and then:
\beq \label{eq:cas1}
\sum_{i=1}^n \arcsin(S_i/S)=\pi~. \eeq
\item Or $x_0$ is not in the interior of $p$, and then:
\beq \label{eq:cas2}
\sum_{i=2}^n \arcsin(S_i/S)=\arcsin(S_1/S)~. \eeq
\end{itemize}
Moreover, in both cases, we have:
$$ \sum_{i=2}^n \arcsin(S_i/S) < \pi~, $$
and:
$$ S_1/S = \sin\left(\sum_{i=2}^n \arcsin(S_i/S)\right)~. $$

Set:
$$ F(s) := s\sin \left( \sum_{i=2}^n \arcsin\left(\frac{S_i}{s}\right)
\right)~. $$
Then:
\begin{eqnarray*}
  F'(s) & = & \sin \left( \sum_{i=2}^n \arcsin\left(\frac{S_i}{s}\right)
\right) + s \cos \left( \sum_{i=2}^n \arcsin\left(\frac{S_i}{s}\right)
\right) \left( \sum_{i=2}^n \frac{-S_i/s^2}{\sqrt{1-S_i^2/s^2}}\right) \\
& = & \cos \left( \sum_{i=2}^n \arcsin\left(\frac{S_i}{s}\right)
\right) \left[ \tan\left( \sum_{i=2}^n
    \arcsin\left(\frac{S_i}{s}\right)\right) - \sum_{i=2}^n \tan\left(
\arcsin\left( \frac{S_i}{s}\right) \right) \right]~.
\end{eqnarray*}
When $\sum_{i=2}^n \arcsin(S_i/s) \in [\pi/2,\pi)$, 
both terms in the first equation are non-positive, and it follows that $F'(s)$
is positive. When $\sum_{i=2}^n \arcsin(S_i/s)\in (0, \pi/2)$, the
second equation shows the same result, because the cosine is positive, while
$\tan$ is convex on $(0,\pi/2)$, so that the second product is also positive. 
Therefore, $F$ is increasing, and it
follows that there is at most one $s$ such that $F(s)=S_1$, and therefore at
most one possible value of $S=\sin(l)$.
 
But, given $l$, the angles $\alpha_i$ are uniquely determined by the $l_i$,
so that there is at most one polygon with vertices on a circle and with edge
lengths equal to the $l_i$.
\epv

We now consider the hyperbolic case. This is a step in the proof of 
Theorem \ref{tm:isoper-hyper}.

\blm \label{lm:unique-hyper}
Let $l_1, \cdots, l_n\in \R_+^*$. There is at most one convex polygon $p\in
\cP^c_H(l)$ which has its vertices either on a circle, on a horocycle, or on a
connected component of the set of points at constant distance from a
geodesic. 
\elm

\bpv
We suppose again that $l_1$ is the largest of the $l_i$. 

Suppose that the vertices of $p$ are on a circle of radius $l$ and center
$x_0$ in $H^2$. Then the edge of $p$ of length $l_i$ is ``seen'' from $x_0$
under an angle $\theta_i$, and the sine law for hyperbolic triangles shows
that: 
$$ \sin(\theta_i/2) = \frac{\sinh(l_i/2)}{\sinh(l)}~. $$
Setting $S_i:=\sinh(l_i/2)$ and $S:=\sinh(l)$, we have again two
possibilities: 
\begin{itemize}
\item Either $x_0$ is in the interior of $p$, and then:
$$
\sum_{i=1}^n \arcsin(S_i/S)=\pi~. $$
\item Or $x_0$ is not in the interior of $p$, and then:
$$
\sum_{i=2}^n \arcsin(S_i/S)=\arcsin(S_1/S)~. $$
\end{itemize}
Moreover, in both cases, we have:
$$ S_1/S = \sin\left(\sum_{i=2}^n \arcsin(S_i/S)\right) $$
and, since $\sin$ is concave on $[0, \pi]$:
$$ S_1/S < \sum_{i=2}^n \sin(\arcsin(S_i/S)) = \sum_{i=2}^n S_i/S~,  $$
so that $S_1<\sum_{i=2}^n S_i$. 

Suppose now that the vertices of $p$ are on a horocycle $h$. An elementary
argument in hyperbolic geometry shows that the distance, along $h$, between
the vertices of an edge of length $l_i$ is equal to
$2\sinh(l_i/2)$. 

In the projective model of $H^2$, the horocycle $h$ appears as an ellipse,
tangent to the boundary at infinity at a point $x_0$. In this model, $p$
appears as a convex polygon, with vertices on this ellipse. 
There is exactly one
edge $e$ of $p$ such that $x_0$ is on one side of $e$, while all the other
edges of $p$ are on the other side. Clearly, $e$ has to be the edge of $p$ of
largest lengths, and, in this case: 
$$ S_1 = \sum_{i=2}^n S_i~. $$

Finally, suppose that all vertices of $p$ are on a connected component $c$ of
the set of points at distance $l$ from a geodesic $g_0$. Let $e_1, \cdots,
e_n$ be the edges of $p$; the sine formula for hyperbolic triangles 
shows that, if $a_i$ is the distance between
the orthogonal projections on $g_0$ of $e_i$, then:
$$ \sinh(a_i/2) = \frac{\sinh(l_i/2)}{\sinh(l)}~. $$

\begin{figure}[h] \label{fg:ellipse}
\centerline{\psfig{figure=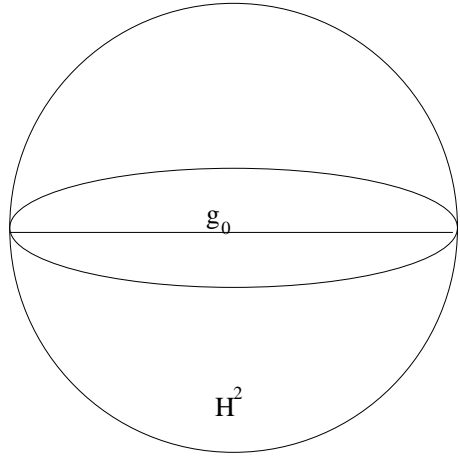,height=6cm}}
\caption{The set of points at constant distance from a geodesic.}
\end{figure} 

\medskip

In the projective model of $H^2$, the set of points at distance $l$ from $g_0$
appears as an ellipse, tangent to the boundary at infinity of $H^2$ at the
endpoints of $g_0$, and $p$ appears as an ellipse (see Figure
\ref{fg:ellipse}). By convexity, there is one
edge $e$ of $p$ such that $g_0$ is on one side, and all the other edges of $p$
on the other. This edge has maximal length, and: 
$$ a_1 = \sum_{i=2}^n a_i~, $$
so that:
$$ \argsinh(S_1/S) = \sum_{i=2}^n \argsinh(S_i/S)~. $$
Then, taking the $\sinh$, we have, because $\sinh$ is convex on $(0,\infty)$:
$$ S_1/S = \sinh \left(\sum_{i=2}^n \argsinh(S_i/S)\right) > \sum_{i=2}^n
  S_i/S~, $$
so that $S_1>\sum_{i=2}^n S_i$.

Considering those three cases, it is clear that the $l_i$ determine whether
there can exist a convex polygon $p$ with edge 
lengths equal to the $l_i$ and vertices
on a circle, a horocycle, or a connected component of the set of points at
constant distance from a geodesic.
\begin{enumerate}
\item If $S_1<\sum_{i=2}^n S_i$, then the only possibility is that the
  vertices of $p$ are on a circle.
\item If $S_1=\sum_{i=2}^nS_i$, then the vertices of $p$ can only be on a
  horocycle. 
\item If $S_1>\sum_{i=2}^n S_i$, the vertices of $p$ can only be on a connected
  component of the set of points at constant distance from a geodesic.
\end{enumerate}
In the first case, we proceed as in the proof of Lemma \ref{lm:unique-sphere} 
above, and set:
$$ F(s) := s\sin \left( \sum_{i=2}^n \arcsin\left(\frac{S_i}{s}\right)
\right)~. $$
Either equation (\ref{eq:cas1}) or equation (\ref{eq:cas2}) is satisfied; in
both cases, $ F(S) = S_1$.
However the same proof as for Lemma \ref{lm:unique-sphere} shows that $F$ is
strictly increasing, so that there is at most one possible value of $S$.
As in the spherical case, it is easy to check that a polygon with vertices on
a circle is uniquely determined (up to global isometries) by its edge lengths
and the radius of the circle; this finishes the proof case (1).

In the second case, the polygon $p$ is clearly uniquely determined (up to
global isometries) by the distances between
its vertices on the horocycle, so the statement in the lemma holds.

In case (3), it is necessary to modify slightly the argument used in the
spherical case. We set: 
$$ G(s) := s\sinh \left( \sum_{i=2}^n \argsinh\left(\frac{S_i}{s}\right)
\right)~. $$
Now:
\begin{eqnarray*}
  G'(s) & = & \sinh \left( \sum_{i=2}^n \argsinh\left(\frac{S_i}{s}\right)
\right) + s \cosh \left( \sum_{i=2}^n \argsinh\left(\frac{S_i}{s}\right)
\right) \left( \sum_{i=2}^n \frac{-S_i/s^2}{\sqrt{1+S_i^2/s^2}}\right) \\
& = & \cosh \left( \sum_{i=2}^n \argsinh\left(\frac{S_i}{s}\right)
\right) \left[ \tanh\left( \sum_{i=2}^n
    \argsinh\left(\frac{S_i}{s}\right)\right) - \sum_{i=2}^n \tanh\left(
\argsinh\left( \frac{S_i}{s}\right) \right) \right]~.
\end{eqnarray*}
But the first factor in the right-hand term is positive, while the second
factor is negative because $\tanh$ is concave on $(0,\infty)$. So $G$ is
strictly decreasing, and there is at most one value of $s$ such that
$G(s)=S_1$. It follows that there is at most one polygon with edge lengths
equal to the $l_i$ and with vertices on a connected component of the set of
points at constant distance from a geodesic.
\epv

\pg{Boundary behavior of the area functional.}

Another, simple element in the proof of Theorem 
\ref{tm:isoper-hyper} is that, on the boundary of the space $\cP^c_S(l)$
(resp. $\cP^c_H(l)$), the interior derivative of the area is positive. 

\blm \label{lm:interior}
Let $p\in \dr \cP^c_S(l)$ (resp. $p\in \dr \cP^c_H(l)$) be a polygon with
exactly one angle equal to $\pi$. Then the interior normal derivative of the
area at $p$ is positive.  
\elm 

\bpv
Let $v$ be the vertex of $p$ at which the angle is equal to $\pi$. Consider
the isometric first-order deformation $U$ of $p$ under which $v$ moves toward
the exterior of $p$, at speed $1$, along the orthogonal to both the edges of
$p$ at $v$. Clearly, $U$ is towards the interior of $\cP^c_S(l)$ (resp. of
$\cP^c_H(l)$) and the first-order variation of the area under $U$ is
positive. The result follows.
\epv

\pg{Proof of Theorem \ref{tm:isoper-hyper}.}

By Lemma \ref{lm:interior}, the area has at least one local maximum in each
connected component of $\cP^c_S(l)$ (resp. of $\cP^c_H(l)$). By Lemma
\ref{lm:critic-hyper}, each critical point of
the area is a polygon with vertices on a circle, a
horocycle, or a connected component of the set of points at constant distance
from a geodesic. By Lemma 
\ref{lm:unique-hyper}, there is at most one such critical point. So the area
has exactly one critical point in $\cP_H^c(l)$, which
is a maximum. By the way, it also follows that $\cP_H^c(l)$ is connected.

Note that the same argument can be used for spherical polygons, it prove the
corresponding result. 

\pg{The convexity of the area.}

In addition, Theorem B$_S$ (resp. B$_H$) shows that the area has a strict
concavity property over at least a subset of the space of convex polygon with
given edge lengths. We call $A$ the function equal to the area of polygons. 

\bdf
Let $p$ be a convex polygon, with vertices $v_1, \cdots, v_n$. Let
$\phi:\R^3\rightarrow \R^3$ be defined by: 
$$ \phi(x) = \sum_{i=1}^n \langle x, v_i\rangle v_i~. $$
Suppose that the $v_i$ are not on a geodesic, let $C_v$ be the vertices of
$p$; we set: 
$$ C_S(p) := \frac{\phi^{-1}(C_v)}{\| \phi^{-1}(C_v)\|}~. $$
\edf

Note that, since the $v_i$ are not on a geodesic, $\phi$ is a self-adjoint,
positive definite operator, so it is invertible. 

To understand the meaning of $C_S(p)$, we consider first the simplest case,
when $p$ is a spherical 
triangle, with vertices $v_1, v_2$ and $v_3$, which we suppose
are not on a geodesic. Then $C_S(p)$ is
the barycenter of the dual polygon $p^*$. Indeed, let 
$p^*=(v_1^*, v_2^*, v_3^*)$; then $(v_1^*, v_2^*, v_3^*)$ is the basis of
$\R^3$ which is dual to the basis $(v_1, v_2, v_3)$. By definition, $C_S(p)$
is defined by the equation:
$$ \frac{\sum_{i=1}^3 \langle v_i, C_S(p)\rangle v_i}{\left\|
  \sum_{i=1}^3 \langle v_i, C_S(p)\rangle v_i \right\| } = C_v(p)~. $$
Taking the scalar product with $v_j^*$, we obtain that $C_S(p)$ is
  characterized by the existence of a $\lambda>0$ such that, for all $j\in \{
  1,2,3\}$:
$$ \langle C_S(p), v_j\rangle = \lambda \langle C_v(p), v_j^*\rangle~. $$
The action of the duality on this equation shows that $C_S(p)$ is the
barycenter of the polygon $p^*=(v_1^*, v_2^*, v_3^*)$.

This example can be extended by taking a triangle with integer ``weights'',
i.e. we consider a polygon with $k_1+k_2+k_3$ vertices, $k_1$ of them being
``collapsed'' on a point $v_1$, $k_2$ on $v_2$, and $k_3$ on $v_3$. The
argument given above carries over to this case, and, again in this case, the
point $C_S(p)$ which is obtained is the barycenter of the triangle dual to
$(v_1, v_2, v_3)$; in other terms, the weights $k_1, k_2, k_3$ have no
influence. This illustrate an interesting ``stability'' property of $C_S(p)$. 

\pg{The concavity of the area.}

Theorem \ref{tm:isoper-hyper} 
is related to the fact that,
over at least part of $\cP_H^c(l)$ (resp. $\cP_S^c(l)$), the area is strictly
concave. 
 
\blm \label{lm:convexity}
Let $p$ be a convex polygon in $S^2$ (resp. $H^2$) with edge lengths given by
$l=(l_1, \cdots, l_n)$. Suppose that $C_S(p)$ is contained in the interior of
$p$. Then the restriction of $A$ to $\cP^c_S(l)$ is strictly 
concave at $p$
for the metric induced by the immersion of $\cP_S(l)$ in $[0,\pi]^n$
defined by sending a polygon to its angles.  
\elm

\bpv
First note that $A$ has an extension as a linear function over $[0,\pi]^n$
since, by the Gauss-Bonnet formula, $A=2\pi-\sum_{i=1}^n(\pi-\alpha_i)$ for
spherical polygons, while $A=-2\pi+\sum_{i=1}^n (\pi-\alpha_i)$ for hyperbolic
polygons. We still call this function $A$, then $dA=(1, \cdots, 1)$ for
spherical polygons, while $dA = -(1, \cdots, 1)$ for hyperbolic polygons.

Let $\hess(A)$ be the Hessian of $A$ over $\Phi_l(\cP_S(l))$. Let $U$ be a
first-order isometric deformation of $p$, then:
$$ \hess(A)(U,U) = dA(\II(d\Phi_l(U), d\Phi_l(U))) 
= \langle \II(d\Phi_l(U), d\Phi_l(U)), I\rangle~, $$
where $\II$ is the second fundamental form of $\Phi_l(\cP^c_S(l))$ and $I$ is
the vector $(1, \cdots, 1)$.
Since $\II(d\Phi_l(U), d\Phi_l(U))$ is orthogonal to $\Phi_l(\cP^c_S(l))$: 
$$ \hess(A)(U,U) = \langle \II(d\Phi_l(U), d\Phi_l(U)), I_N\rangle~, $$
where $I_N$ is the orthogonal projection of $I$ on the normal space of
$\Phi_l(\cP^c_S(l))$ at $\Phi_l(p)$. By definition of $\Psi_p$,
$I_N=\Psi_p(w)$, where $w\in \R^3$ is characterized by the fact that
$\Psi_p(w)=I_N$, i.e. that $\Psi_p(w)-I\in T_{\Phi_l(p)}\Phi_l(\cP^c_S(l))$,
which translates as:
$$ \sum_{i=1}^n (\langle w, v_i\rangle -1)v_i = 0~, $$
so that:
$$ \sum_{i=1}^n \langle w, v_i\rangle v_i = \sum_{i=1}^n v_i~. $$
This shows that $w= \lambda C_S(p)$, for some $\lambda >0$. 

Applying Remark \ref{rk:II} then leads to:
$$ \hess(A)(U,U) =
- \lambda \langle b(U), C_S(p)\rangle~, $$
so that the Hessian of $A$ is negative definite at $p$. 

The proof in the spherical case is almost the same, with some differences in
the signs, and we leave it to the interested reader.
\epv

Although the proof given here only works when $C_S(p)$ is contained in the
interior of $p$, it is not obvious whether the result is valid for all convex
polygons -- the existence of a unique critical point of the area over
$\cP_S^c(l)$, which is true in all cases, suggest that it might be the case.

\section{Infinitesimal rigidity of Euclidean polyhedra}

\pg{Infinitesimal rigidity of convex polyhedra.}

In this section we consider convex polyhedra in $\R^3$. We define a convex
polyhedron as the boundary of a compact subset of $\R^3$ which is the
intersection of a finite number of half-spaces. Given a polyhedron $P$, which
is the boundary of a compact subset $\Omega\subset \R^3$, where $\Omega$ is
the intersection of distincts closed 
half-spaces $H_1, \cdots, H_n$ (each of which intersect the interior of
$\Omega$, a vertex of $P$
is a point of $P$ which is contained in at least 3 of the boundaries of the
$H_i$. An edge of $P$ is a maximal connected 
subset of $P$ which is contained in the boundaries of two of the
$H_i$ but containing no vertex. A face of $P$ is a connected component of the
complement of the vertices and the edges.

The geometry of Euclidean polyhedra has interested geometers for quite a long
time \cite{euclid}.
A nice result of Legendre \cite{legendre}\footnote{The contribution of
  Legendre was recently pointed out by I. Sabitov} and Cauchy \cite{cauchy}
states that convex polyhedra are rigid: if two convex polyhedra have the same
combinatorics and the same induced metric on their faces, 
they are congruent. A related
result, first proved by Dehn \cite{dehn-konvexer} 
(which is also a consequence of the ideas
of Legendre and Cauchy) is that convex polyhedra are infinitesimally rigid:
any first-order deformation of a polyhedron $P$ which acts as a Killing field
on each of its faces is the restriction to $P$ of a global Killing field.

From a modern viewpoint, one of the main motivation to study the infinitesimal
rigidity of convex polyhedra is that it is the key point in the proof of the
following nice theorem. 

\begin{thm}[Aleksandrov \cite{alex}]
The induced metric on any convex polyhedron in $\R^3$ is flat, with conical
singularities where the total angle is less than $2\pi$. Conversely, any flat
metric on $S^2$ with conical singularities where the total angle is less than
$2\pi$ is induced on a unique (up to isometries) convex polyhedron in $\R^3$.
\end{thm}

To prove this theorem, one considers the natural map sending a convex
polyhedron to its induced metrics; it is a map between two spaces of the same
dimension, and the infinitesimal rigidity statement means that its
differential is everywhere an isomorphism, so that it is a local
homeomorphism.

We will show here that Theorem C is a consequence
of Theorem B$_S$, thus providing a new proof which is rather natural. It will
also serve as an introduction for the rigidity result of the next section,
since the proof given there is related but more complicated. The rigidity
proof given here 
is related to a proof of the Cauchy-Legendre rigidity result discovered by
Pogorelov \cite{pogo-polygones}.

\pg{Closed convex polyhedra.} 

Let $P$ be a closed, convex Euclidean polyhedron. Choose a point $p_0$ in the
interior of $P$. Let $u_0$ be the function defined on $\R^3$ by:
$u_0(p)=d(p_0,p)^2/2$. We consider a first-order deformation of $P$, i.e. a
vector field $v$ on $P$ which acts as a Killing field on each of its faces; we
will show that $v$ is trivial, i.e. it acts as a Killing field on $P$ as a
whole. 

Let $e$ be an edge of $P$. The first-order deformation $v$ of $P$ determines a
first-order variation $\thetad_e$ of the dihedral angle $\theta_e$ of $P$ at
$e$. $v$ also determines a function, which we call $\ud_0$, on $P$, defined as
the first-order variation of the restriction of $u_0$ to $P$. By definition of
$u_0$, we have:
\beq \label{eq:u0}
\forall p\in P~,~~ \ud_0(p) = \langle v(p),p_0p\rangle~. \eeq 
We define a 1-form $W$ on the 1-skeleton of $P$ (the union of its
edges) as follows: for each vector $u$ tangent to an edge $e$, $W(u)=\thetad_e
d\ud_0(u)$. 

Consider a parametrization of $e$ at constant speed, say: $p(t)=at+b$, with
$a,b\in \R^3$. Then, using (\ref{eq:u0}):
$$ \ud_0(p(t))' = d\ud_0(p'(t)) 
= \langle dv(p'(t)), p_0p(t)\rangle + \langle v(p(t)), p'(t)\rangle~. $$
But $p'(t)=a$ is independent of $t$, and $dv(p'(t))$ does not depend on $t$
  either since $v$ acts on $e$ as a Killing field. So:
$$ \frac{d}{dt}  \ud_0(p(t))' = 2 \langle dv(p'(t)), p'(t)\rangle~. $$
But $v$ acts isometrically on $e$, so that $\langle dv(p'(t)), p'(t)\rangle
=0$, and it follows that $d\ud_0(p'(t))$ is constant over $e$, so $\thetad_e
d\ud_0(p'(t))$ is constant over $e$.

Now consider an oriented edge $e$ of $P$, we call $e_-$ and $e_+$ its
endpoints, and we let $W_e$ be the number $W(p'(t)=\thetad_e
d\ud_0(p'(t))$, where $p(t)$ is a parametrization of $e$ at speed $1$,
respecting the orientation. Then we have:
\beq \label{eq:sum0} \sum_{v} \sum_{e_-=v} W_e = 0~, \eeq
because each non-oriented edge contributes twice, with opposite signs. 

Let $x$ be one of the vertices of $P$. Let $e_1, \cdots, e_{n-1},
e_n=e_0$ be the oriented edges of $P$ such that $e_{i,-}=x$, 
with their natural cyclic
order, and let $v_1, \cdots, v_n=v_0$ be unit vectors with $v_i$ in the
direction of $e_i$, oriented from $x$ towards the other vertex of
$e_i$. Then $v_1, \cdots, v_n$ are the vertices of a spherical 
polygon (called the link
of $P$ at $x$) which is convex since $P$ is convex. 

For each $i\in \{ 1,\cdots, n\}$, the first-order deformation $v$ of $P$
determines a first-order deformation $v'$ of the spherical 
polygon $(v_1, \cdots, v_n)$, given
at each vertex $v_i$ by the vector $\vd_i = dv(v_i)$; note that $dv(v_i)$ is
orthogonal to $v_i$ because $v$ is an isometric deformation. Since $v$ is
isometric, it does not change (at first order) the interior angles of the
faces of $P$; since the interiors angles of the faces of $P$ at $x$ are the
lengths of the edges of $(v_1, \cdots, v_n)$, the first-order deformation $v'$
is isometric. 

The angles $\alpha_1, \cdots, \alpha_n$ of the polygon $(v_1, \cdots, v_n)$ are
equal to the dihedral angles $\theta(e_1), \cdots, \theta(e_n)$ of $P$ at the
edges $e_1, \cdots, e_n$. Therefore: 
\begin{eqnarray*}
  \sum_{i=1}^n W_{e_i} & = & \sum_{i=1}^n \left( \langle dv(v_i), p_0x\rangle +
  \langle v(x), v_i\rangle \right) \thetad_i \\
& = & \sum_{i=1}^n \left( \langle \vd_i, p_0x\rangle +
  \langle v(x), v_i\rangle \right) \alphad_i \\
& = & \left\langle  \sum_ {i=1}^n \alphad_i \vd_i, p_0x\right\rangle +
  \left\langle v(x), 
  \sum_{i=1}^n \alphad_i v_i\right\rangle~. 
\end{eqnarray*}
Since $p_0$ is in $P$ and $P$ is convex, $-p_0x$ is in the positive cone over
the interior of $(v_1,
\cdots, v_n)$, so that the first sum is non-positive by Theorem B$_S$, and it
vanishes if and only if the deformation $v'$ is trivial. The second sum is
zero by Theorem A$_S$. So, by (\ref{eq:sum0}), the first sum is $0$ for each
vertex $x$ of $P$, and none of the angles of $P$ varies in the first-order
deformation $v$, so that $v$ is a trivial deformation. 

\pg{Further results ?}

The proof of Theorem C given here has some flexibility. For
instance, it should be possible to consider polyhedral surfaces with some
singularities, e.g. ramifications points inside some of the faces. One
can ask whether it also applies for polyhedral surfaces having a kind of
ramification at the vertices, or to polyhedral surfaces which are not convex
but satisfy a kind of weak convexity like the one considered in
\cite{stoker,kann,rodriguez-rosenberg}. 

Another possible extension would be to convex surfaces with boundary, as
considered e.g. in \cite{kann}; in this setting it is necessary to have some
constraints on the deformations at the boundary vertices. However, only the
local convexity at the vertices plays a role, so that it should be possible to
consider surfaces of higher Euler characteristic. 

\section{Convex polyhedral surfaces in the Minkowski space}

This section contains an extension of Theorem C to a rigidity question for
convex, equivariant surfaces in the Minkowski space, for which the techniques
used here are particularly well adapted since there is a distinguished
point. We first define the polyhedral surfaces that are considered. Then we
basically follow the path taken above in the Euclidean case, adapting the
proof from the sphere 
to the de Sitter plane. Along the way we recall some elementary facts of
de Sitter geometry, proved here for completeness. 

\pg{Equivariant embeddings of surfaces.}

Let $S$ be a closed surface, which will be of genus at least $2$ in all this
section. We first define equivariant polyhedral embeddings of $S$. More general
definitions could of course be given, but we stick to what is really necessary
in our context.

\bdf
An {\it equivariant space-like polyhedral embedding} of $S$ in $\R^3_1$ 
is a couple $(\phi, \rho)$, where:
\begin{itemize}
\item $\phi$ is a space-like polyhedral embedding of the universal cover $\St$
  of $S$ in $\R^3_1$, i.e.: 
  \begin{itemize}
  \item $\phi$ is continous. 
  \item There exists a cellular decomposition of $S$ as the
  union of a finite number of cells, each the image by a diffeomorphism of the
  interior of a convex polygon in $\R^2$, 
  such that the image by $\phi$ of each cell
  of the corresponding cellular decomposition of $\St$ is the interior of a
  convex polygon is some space-like plane in $\R^3_1$.
  \item For each space-like plane $H$ in $\R^3_1$, the orthogonal projection
  of $\phi(\St)$ on $H$ is one-to-one. 
  \end{itemize}
\item $\rho$ is a morphism from $\pi_1(S)$ to the isometry group of $\R^3_1$.
\item For each $x\in \St$ and each $\gamma\in \pi_1(S)$, $\phi(\gamma
  x)=\rho(\gamma)\phi(x)$. 
\end{itemize}
We will say that $(\phi,\rho)$ is {\it Fuchsian} if the image of $\rho$ is
contained in the identity component of the sub-group of isometries fixing the
origin, i.e. $\rho(\pi_1(S))\subset SO_0(2,1)$.
\edf

A direct consequence of the definition of an equivariant embedding of $S$ in
$\R^3_1$ is that the metric induced on $\St$ by $\phi$ is invariant under the
action of $\pi_1S$, so that an equivariant embedding of $S$ induces a metric
on $S$.

Consider an equivariant embedding $(\phi,\rho)$ of $S$ in $\R^3_1$. There is a
natural notion of first-order deformation of $(\phi,\rho)$; it corresponds to
the deformations of $\phi$ among the equivariant embeddings of
$S$. Restricting our attention to the deformations among Fuchsian equivariant
embeddings, such a
deformation can be described as a couple $(\phid, \rhod)$, where:
\begin{itemize}
\item $\phid$ is a vector field defined over $\phi(\St)$.
\item $\rhod:\pi_1S\rightarrow TSO(2,1)$ is a map, such that:
$$ \forall \gamma \in \pi_1S, \rhod(\gamma)\in T_{\rho(\gamma)}SO_0(2,1)~, $$
and, for all $\gamma, \gamma'\in \pi_1S$, if we consider $\rho(\gamma)$ and
  $\rho(\gamma')$ as acting on $SO(2,1)$ and $\rhod$ as sending elements of
  $\pi_1S$ to vector fields on $SO(2,1)$, we have:
$$ \rhod(\gamma \gamma') = \rhod(\gamma)(\rho(\gamma')) +
\rho(\gamma)_*\rhod(\gamma')~. $$ 
\item If we identify the elements of $so(2,1)$ with the Killing fields on
  $\R^3_1$ which vanish at $0$, then:
$$ \forall x\in \St, \forall \gamma\in \pi_1S, \phid(\gamma x) =
\rho(\gamma)_* \phid(x) + \rhod(\gamma)\phi(x)~. $$
\end{itemize}

\pg{Fuchsian polyhedral surfaces.}

The main result of this section is the following infinitesimal rigidity
statement. There is an analogous statement for smooth surfaces in \cite{iie}. 

\begin{thm} \label{tm:minkowski}
  Let $(\phi,\rho)$ be a Fuchsian, convex equivariant embedding of a surface
  $S$ (of genus at least $2$) in $\R^3_1$. Let $(\phid, \rhod)$ 
  be a first-order
  deformation of $(\phi, \rho)$ among Fuchsian equivariant embeddings. If the
  first-order variation of the metric induced on $S$ by $(\phi,\rho)$
  vanishes, then $(\phid, \rhod)$ is trivial, 
  i.e. $\phid$ is the restriction to $\phi(\St)$ of
  a Killing field of $\R^3_1$ which vanishes at $0$, and $\rhod=0$.
\end{thm}

The proof comes after some preliminaries, basically following the proof of
Theorem C and checking that the various parts of the proof carry over from
the Euclidean to the Minkowski 3-dimensional space, and from the sphere to the
de Sitter plane. 

As in the Euclidean case, the main motivation for this theorem is that it
should be a key point in the proof a statement describing the metrics induced
on fuchsian equivariant polyhedra in $\R^3_1$, as in the Aleksandrov theorem
quoted above for the Euclidean space. Namely, one would like to answer the
following question: 
\begin{quote} {\it
  Let $S$ be a surface of genus at least $2$, and let $g$ be a flat metric on
  $S$ with conical singularities, with total angle larger than $2\pi$ at each
  singular point. Is there a unique (up to global isometries) convex
  Fuchsian polyhedral embedding of $S$ in $\R^3_1$ such that the induced
  metric is $g$ ? }
\end{quote}
We do not answer this question here since it would demand some considerations
leading us too far from Theorems A$_S$ and B$_S$, the main theme of this text. 

\pg{Distances and angles in the de Sitter plane.}

As mentioned in the introduction, the de Sitter plane, denoted by $S^2_1$,
is the quadric of
equation $\langle x,x\rangle =1$ in $\R^3_1$, with the induced metric. Details
can be found e.g. in \cite{coxeter-de-sitter,coxeter-non-euclidean,O}.

Let $x,x'$ be two points in $S^2_1$. By analogy with the distance in
the sphere $S^2\subset \R^3$, 
we define the ``distance'' between $x$ and $x'$ as the number $d(x,y)\in
\C/2\pi \Z$ such that:
$$ \langle x,x'\rangle = \cos(d(x,x'))~. $$
This leaves an indetermination concerning the sign of $d(x,x')$, which is
removed by the following explicit description:
\begin{itemize}
\item $d(x,x')\in (0,\pi)$ when $x$ and $x'$ are distinct, non-antipodal
  points on a space-like geodesic.
\item $d(x,x')=0$ when $x$ and $x'$ are on a light-like geodesic.
\item $d(x,x')\in i\R_+$ when $x$ and $x'$ are on a time-like geodesic.
\item $d(x,x')\in \pi-i\R_+$ when $x$ and $x'$ are in two different connected
  components of the intersection of $S^2_1$ with a time-like plane containing
  $0$ (each connected component of the intersection 
  is then a time-like geodesic).
\end{itemize}
This definition differs from the one used in \cite{shu} by a factor $i$,
basically because here we stick as close as possible to spherical geometry,
while the emphasis in \cite{shu} was on hyperbolic geometry. The definition of
\cite{shu} would be obtained if we had taken the $\cosh$ of $d(x,x')$ instead
of the cosine in the definition above. The proof below would then have less
factors ``$i$'', but it would be based on the hyperbolic rather
than the spherical trigonometric 
formulas, thus blurring the analogy with the proof of Theorem C.

There is a related notion of angles in $S^2_1$. Let $x\in S^2_1$, and let
$u,v\in T_xS^2_1$ be such that $\langle u,u\rangle\neq 0$ and $\langle
v,v\rangle \neq 0$; the angle between them is a number $\theta\in \C/2\pi\Z$
such that:
$$ \langle u,v\rangle^2 = \cos^2(\theta) \langle u,u\rangle \langle
  v,v\rangle~. $$
Note that $\theta$ is defined up to sign and up to the transformation
  $\theta\mapsto \pi-\theta$, which is quite normal since it depends on the
  direction of $u$ and $v$ and on whether we consider $u$ or $v$ first.
This ambiguity which is removed by the
following explicit description, in which we suppose that $(u,v)$ is a
positively oriented basis of $T_xS^2_1$:
\begin{itemize}
\item $\theta\in i\R_+$ if $u$ and $v$ are both future-oriented time-like
  vector (or if they are both past-oriented time-like vectors).
\item $\theta \in \pi - i\R_+$ if $u$ and $v$ are time-like vectors, but $u$ is
  future-oriented and $v$ is past-oriented.
\item $\theta\in i\R_+$ if $u$ and $v$ are both space-like vectors, and are in
  the same connected component of the set of space-like vectors at $x$.
\item $\theta \in \pi - i\R_+$ 
  if $u$ and $v$ are space-like vectors in different
  connected components of the set of space-like vectors at $x$.
\item $\theta \in \pi/2 + i\R$ if $u$ is space-like and $v$ is time-like, or
  conversely.
\end{itemize}

\pg{Triangle formulas in the de Sitter plane.}

With those definitions, equation (\ref{eq:cos}) and the sine formula hold in
$S^2_1$ exactly as in the sphere.

\begin{prop} \label{pr:triangles}
Let $(A,B,C)$ be a triangle in the de Sitter plane, with edge lengths
$a,b,c$ and angles $\alpha, \beta,\gamma$. Suppose that $a, b, c\neq 0$. Then:
$$ \cos(a) = \cos(b)\cos(c) + \cos(\alpha)\sin(b)\sin(c)~. $$
Moreover, if $a,b,c\neq 0$, then:
$$ \frac{\sin(\alpha)}{\sin(a)} = \frac{\sin(\beta)}{\sin(b)}
= \frac{\sin(\gamma)}{\sin(c)}~. $$
\end{prop}

\begin{proof}
Recall the definition of the cross-product $X\boxtimes Y$, in coordinates, 
used in the Minkowski space:
$$ (x_1,x_2,x_3)\boxtimes (y_1,y_2,y_3) = (x_2y_3-x_3y_2, x_3y_1-x_1y_3,
-x_1y_2+x_2y_1)~. $$
Note that this definition is natural insofar as it can be obtained like the
vector product in Euclidean space, by associating 1-forms to vectors, taking
the wedge product, and associating a vector to the resulting 2-form. This
shows that the definition given using coordinates is independent of the
orthonormal basis of $\R^3_1$ which has been used.

It follows from this definition that, given 3 vectors $X,Y$ and $Z$ in $\R^3$,
the number $\langle X,Y\boxtimes Z\rangle$ is the same whether the cross
product 
and the scalar product are considered for the Euclidean or the Minkowski space
structure (there are two sign differences but they cancel out). 
Therefore, it remains true in the Minkowski space that:
\beq \label{eq:triple}
\langle X,Y\boxtimes Z\rangle = \langle Y, Z\boxtimes X\rangle = \langle Z, X\boxtimes
Y\rangle~. \eeq
The same argument also shows that $\langle X, X\boxtimes Y\rangle = \langle Y,
X\boxtimes Y\rangle =0$ in the Minkowski space, so that $X\boxtimes Y$ is orthogonal
(for the Minkowski scalar product) to $X$ and to $Y$.

It is also easy to check, by taking a ``good'' orthonormal basis of $\R^3_1$,
that, if $X$ is a space-like unit vector orthogonal to $Y$ and $Z$, then:
$$ \langle X\boxtimes Y, X\boxtimes Z\rangle = \langle Y,Z\rangle~. $$
Moreover, let $X,Y\in S^2_1$ be two non-collinear vectors of $\R^3_1$, 
there exists a direct orthonormal basis $(X,Y',Z)$ of $\R^3_1$, such that
$Y=\cos(d(X,Y))X + \sin(d(X,Y))Y'$; a direct computation then shows that:
$$ X\boxtimes Y = \sin(d(X,Y)) Z~. $$

Now let $A,B,C\subset S^2_1$ be the 3 vertices of a de Sitter triangle. 
Then $B$ and $C$ can be
decomposed as $B=B_\parallel A + B_\perp$ and $C=C_\parallel A + C_\perp$,
where $B_\parallel, C_\parallel\in \R$ and $B_\perp, C_\perp\in \R^3_1$ 
are vectors orthogonal to $A$. Then: 
$$ \langle A\boxtimes B, A\boxtimes C\rangle = \langle A\boxtimes (B_\parallel
A + 
B_\perp), A\boxtimes (C_\parallel A + C_\perp) \rangle = \langle A\boxtimes
B_\perp, 
A\boxtimes C_\perp \rangle = \coupeq  
= \langle B_\perp, C_\perp\rangle = \langle
B,C\rangle - B_\parallel C_\parallel \langle A,A\rangle = \langle
B,C\rangle - \langle A,B\rangle \langle A,C\rangle~. $$
But $A\boxtimes B=\sin(c) N_C$ and $A\boxtimes C = \sin(b)N_B$, where $N_C$ and
$N_B$ are the unit vectors orthogonal to the (oriented) plane containing
$0,A,B$ and $0,A,C$, respectively. By definition, $\langle N_C,N_B\rangle =
\cos(\alpha)$, and we find that:
$$ \cos(\alpha)\sin(b)\sin(c) = \langle A\boxtimes B, A\boxtimes C\rangle =
\langle B,C\rangle - \langle A,B\rangle \langle A,C\rangle = \cos(a) -
\cos(b)\cos(c)~, $$
which proves the first equation.

To prove the sine formula, 
note that $N_B$ and $N_C$ are both orthogonal to $A$,
so that $(A\boxtimes B)\boxtimes (A\boxtimes C)$ is collinear to $A$. It
follows that:  
$$ \langle (A\boxtimes B)\boxtimes (A\boxtimes C), A\rangle = \sin(c)
\sin(b)\sin(\alpha)~. $$
But, using the same decomposition of $B$ and $C$ as above, we have:
$$ (A\boxtimes B)\boxtimes (A\boxtimes C) = (A\boxtimes B_\perp)\boxtimes
(A\boxtimes C_\perp) = 
B_\perp \boxtimes C_\perp = B\boxtimes C - B_\parallel (A\boxtimes C_\perp) -
C_\parallel (A \boxtimes B_\perp)~; $$
taking the scalar product with $A$ yields:
$$ \langle (A\boxtimes B)\boxtimes (A\boxtimes C), A\rangle = \langle
B\boxtimes C,
A\rangle~. $$
It follows that the quantity $\sin(b)\sin(c)\sin(\alpha)$ is invariant under
a cyclic permutation on $(a,b,c)$ and $(\alpha, \beta, \gamma)$, and the
second equation in the proposition follows. 
\end{proof}

\pg{Polygons in the de Sitter plane.}

The analog of Theorem A$_S$ and Theorem A$_H$ also holds for de Sitter
polygons with non-degenerate edges (when one of the edges is degenerate, the
notion of angle is not well-defined, so that the statement would not make
sense). 

\begin{thn}{A$_{dS}$}
  Let $p=(v_1, \cdots, v_n)$ be a de Sitter polygon with non-degenerate
  edges. Let $\alpha_1, \cdots,
  \alpha_n$ be its angles, and let $\alphad_1, \cdots, \alphad_n\in \R$
  be a first-order variation of its angles induced by an isometric first-order
  deformation of $p$. Then:
$$ \sum_{i=1}^n \alphad_i v_i=0~, $$
where the $v_i$ are considered as points in $S^2_1\subset \R^3_1$. 

Conversely, if this
equation is satisfied by an $n$-uple $(\alphad_1, \cdots,\alphad_n)$ and
moreover the $v_i$ are not all on a de Sitter geodesic -- or on
the intersection of
$S^2_1$ with a Minkowski plane containing the origin -- then there exists an
isometric first-order deformation of $p$ such that the $\alphad_i$ are the
associated first-order variations of the $\alpha_i$.
\end{thn}

The proof follows exactly the proof of Theorem A$_H$.

Again as in the sphere, Theorem A$_{dS}$ leads to the definition of a
quadratic invariant $b$ defined on the first-order isometric deformations of a
de Sitter polygon (with non-degenerate edges). $b$ is defined as:
$$ b(U) = \sum_{i=1}^n d\alpha_i(U)dv_i(U)~, $$
As in the spherical and the hyperbolic case, if two first-order
infinitesimal deformations $U$ and $U'$ differ by a trivial deformation, then
$b(U)=b(U')$. However, the definition of the angles which we have used means
that $b(U)$ is imaginary. 

\pg{The positivity of $b$.}

In the same way, the ``positivity'' property of $b$ which we found in the
spherical case still holds for the de Sitter polygons which are duals of
convex hyperbolic polygons. 

\begin{thn}{B$_{dS}$}
Let $p$ be a convex polygon in $S^2_1$, which is dual to a convex hyperbolic
polygon $p^*$, and let $U$ be a non-trivial infinitesimal
first-order deformation of $p$. Then $ib(U)\in \R_+^*\mbox{int}
(p^*)$, i.e. $ib(U)$ is contained in the positive cone over the interior of 
$p^*$.
\end{thn}

The proof follows the proof of Theorem B$_S$, but some additional details are
necessary. Using Proposition \ref{pr:triangles}, the computations done in
section 3 carry over to the de Sitter case. In particular, the scalar product
of $b(U)$ with $v_1$ remains diagonal in the basis $(U_2, \cdots, U_{n-2})$. 
Equation (\ref{eq:calcul}) can be
rewritten in a way which is more convenient for us, as:
$$ \langle b(U_i), v_1\rangle 
= \frac{\sin^2(d_{1,i+1})\sin(\alpha'_1)}{\sin(d_{0,i+1}) \sin(d_{i,i+1})
  \sin(\alpha'_0)\sin(\alpha'_i)} (dd_{1,i+1}(U_i))^2~. $$
This follows from (\ref{eq:calcul}) because, in Lemma \ref{lm:calcul}, the
definition of $V$ is that, under this first-order variation, the distance
$d_{1,i+1}$ 
between $v_1$ and $v_{i+1}$ varies at speed $1$; this also uses the fact that
$b$ is a quadratic form. Recall that, in this equation $\alpha'_0, \alpha'_1$
and $\alpha'_i$ are not the angles of the polygon $p$, but rather the angles
of the quadrilateral $(v_0,v_1, v_i, v_{i+1})$ at the vertices $v_0, v_1$ and
$v_i$. 

Using as variable the $\cos$ of the distance from
$v_1$ and $v_{i+1}$, we obtain: 
\beq \label{eq:calcul2}
\langle b(U_i), v_1\rangle = \frac{\sin(\alpha'_1)}{\sin(d_{0,i+1}) 
\sin(d_{i,i+1}) \sin(\alpha'_0)\sin(\alpha'_i)} (d\langle v_1, v_{i+1}\rangle
(U_i))^2~. \eeq

Clearly, $(d\langle v_1, v_{i+1}\rangle (U_i))^2\in \R_+$; moreover, since the
edges of $p$ are space-like, $d_{i,i+1}\in (0,\pi)$, so that
$\sin(d_{i,i+1})\in (0,1)$. 
However, to understand the signs of the various sines appearing in this
equation, it is necessary to consider four different cases, depending on
whether $d_{1,i}$ and $d_{0,i+1}$ are in $(0,\pi)$ or in $\pi-i\R_+^*$, 
and then the limit cases. 
The four cases are shown, in the projective model of half of
$S^2_1$, in Figure 2.  %%% \ref{fg:projectif}.

\begin{figure}[h] \label{fg:projectif}
\centerline{\psfig{figure=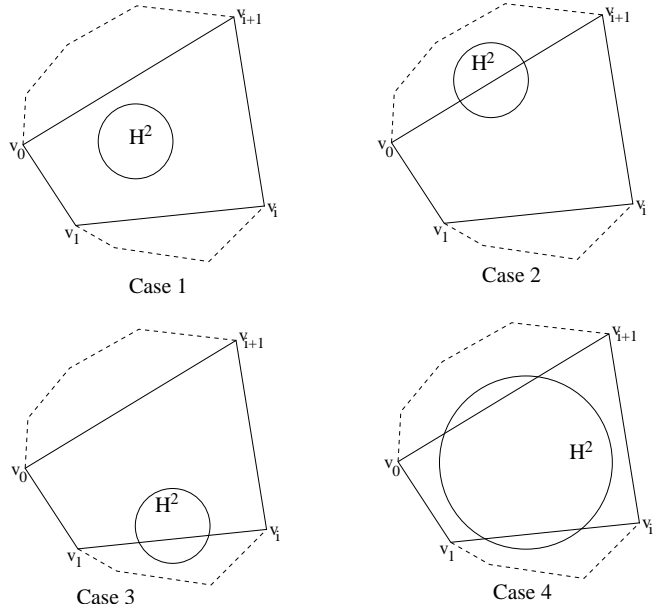,height=8cm}}
\caption{Four possible cases.}
\end{figure} 

\medskip

\noindent{\bf 1st case: $d_{1,i}\in (0,\pi), d_{0,i+1}\in (0,\pi)$.}
In other terms, $v_0$ and $v_{i+1}$ are on a space-like geodesic, and $v_1$ and
$v_i$ are on a space-like geodesic.
Then $\sin(d_{0,i+1})\in (0,1)$. Moreover, the angles $\alpha'_0,\alpha'_1$
and $\alpha'_i$ are all of the form $\pi-ir$, for some $r>0$, so that
$\sin(\alpha'_0), \sin(\alpha'_1)$ and $\sin(\alpha'_i)$ are in
$i\R_+^*$. So equation (\ref{eq:calcul2}) shows that $\langle
b(U_i), v_1\rangle \in -i\R_+^*$.

\medskip

\noindent{\bf 2nd case: $d_{1,i}\in (0,\pi), d_{0,i+1}\in \pi-i\R_+^*$.} 
This means that $v_1$ and $v_i$ are on a space-like geodesic,
  while $v_0$ and $v_{i+1}$ are on the intersection with $S^2_1$ of a
  time-like plane in $\R^3_1$ containing $0$ (but on different connected
  components of this intersection).
  Now $\sin(d_{0,i+1})\in i\R_+^*$. $\alpha'_1$ and
  $\alpha'_i$ are of the form $\pi-ir$ for some $r>0$, 
  so that $\sin(\alpha'_1)$ and
  $\sin(\alpha'_i)$ are in $i\R_+^*$, 
  but $\alpha'_0$ is of the form $\pi/2
  +ir$, for some $r\in \R$, and it follows that $\sin(\alpha'_i) =
  \cos(ir)=\cosh(r) \in
  [1,\infty)$. Therefore, (\ref{eq:calcul2}) again shows that $\langle
  b(U_i), v_1\rangle \in -i\R_+^*$.

\medskip

\noindent{\bf 3rd case: $d_{1,i}\in \pi-i\R_+^*, d_{0,i+1}\in (0,\pi)$.} 
Then $\sin(d_{0,i+1})\in (0,1)$. Moreover,
  $\alpha'_0$ is of the form $\pi-ir$, for some $r>0$, while $\alpha'_1$ and
  $\alpha'_i$ are of the form $\pi/2+ir$, for some $r\in \R$. So
  $\sin(\alpha'_0)\in i\R_+^*$, while
  $\sin(\alpha'_1),\sin(\alpha'_i)\in [1,\infty)$. By (\ref{eq:calcul2}),
  $\langle b(U_i), v_1\rangle \in -i\R_+^*$ also in this case.

\medskip

\noindent{\bf 4th case: $d_{1,i}\in \pi - i\R_+^*, d_{0,i+1}\in
  \pi-i\R_+^*$.}  
In this case, $\sin(d_{0,i+1})\in i\R_+^*$, and
$\alpha'_0, \alpha'_1$ and $\alpha'_i$ are of the form $\pi/2+ir$, for some
$r\in \R$. Therefore, $\sin(\alpha'_0),\sin(\alpha'_1)$ and $\sin(\alpha'_i)$
are in $[1,\infty)$, and, again,  $\langle b(U_i), v_1\rangle \in
-i\R_+^*$.

\medskip

\noindent{\bf Last case: $v_1$ and $v_i$ are in the intersection with $S^2_1$
  of a light-like plane in $S^2_1$ containing $0$, or the same holds of $v_0$
  and $v_{i+1}$.} The result the follows from an approximation of the
  quadrilateral $(v_0,v_1,v_i,v_{i+1})$ by a sequence of quadrilaterals which
  are in one of the four cases detailed above, and a corresponding sequence of
  approximations of $U_i$.

\medskip

Summing over $i\in \{ 2, \cdots, n-2\}$, we find that:
$$ \langle ib(U), v_1\rangle = \sum_{i=2}^{n-2} i\langle b(U_i), v_1\rangle
>0~. $$ 
The same holds with $v_1$ replaced by any of the other vertices of $p$, and
Theorem B$_{dS}$ follows. 

\pg{Proof of the rigidity theorem.}

The proof of Theorem \ref{tm:minkowski} follows quite precisely the proof of
Theorem C. We consider a Fuchsian, equivariant, convex polyhedral embedding
$(\phi, \rho)$ of
a surface $S$ of genus at least $2$ in $\R^3_1$. Let $p_0=0\subset \R^3_1$,
i.e. $p_0$ is the point which is fixed by $\rho$. We call $u_0$ the function,
defined on $\R^3_1$, as: $u_0(x) = \langle x,x\rangle/2$. Then, since $0$ is
fixed by the representation $\rho$, the restriction of $u_0$ to $\phi(\St)$ is
invariant under the action of $\pi_1S$ by $\rho$, so that $u_0$ defines a
function over $S$.

Let $(\phid, \rhod)$ be a first-order deformation of $(\phi,\rho)$ among
Fuchsian embeddings. Let $e$ be an oriented edge of $\phi(\St)$, let
$\theta_e$ be the corresponding dihedral angle, and let be $\thetad_e$ be the
first-order variation of this dihedral angle under the first-order variation
$(\phid,\rhod)$. We associate to $e$ the number $W_e$ defined as: $W_e =
\thetad_e d\ud_0(p'(t))$, where $p(t)$ is a parametrization of $e$ at speed
one. As in the Euclidean case -- and for the same reasons -- $W_e$ is
independent of $t$. By construction, this quantity is invariant under the
action of $\pi_1S$. 

Again in this setting, $\sum_v \sum_{e_-=v} W_e=0$. On the other hand, let $x$
be a vertex of $\phi(\St)$. Let $e_1, \cdots, e_n$ be the oriented edges
starting from $x$, in the cyclic order in which they appear. The link of
$\phi(\St)$ at $x$ is a convex space-like polygon $p$, 
with vertices $v_1, \cdots, v_n$, corresponding to $e_1, \cdots, e_n$, 
which is the dual  of a hyperbolic polygon. 
The first-order isometric deformation $(\phid, \rhod)$ induces a first-order
deformation $U$ of $p$. Since $(\phid, \rhod)$ is
isometric, it does not change (at first order) the interior angles of the
faces of $\phi(\St)$ adjacent to $x$, so that the first-order deformation
$U$ is isometric. 

As in the Euclidean case, if $\vd_i$ is the first-order displacement of $v_i$
under $U$ and if $\alpha_i$ is the angle of $p$ at $v_i$, then:
\begin{eqnarray*}
  \sum_{i=1}^n W_{e_i} & = & \sum_{i=1}^n \left( \langle d\phid(v_i), 
x\rangle +
  \langle \phid(x), v_i\rangle \right) \thetad_{e_i} \\
& = & \sum_{i=1}^n \left( \langle \vd_i, x\rangle +
  \langle \phid(x), v_i\rangle \right) \alphad_i \\
& = & \left\langle  \sum_ {i=1}^n \alphad_i \vd_i, x\right\rangle +
  \left\langle \phid(x), 
  \sum_{i=1}^n \alphad_i v_i\right\rangle~. 
\end{eqnarray*}
The second terme vanishes by Theorem A$_{dS}$. 
By convexity of $\phi(\St)$, $x$ is in the interior of $p$, while, 
by Theorem B$_{dS}$, $ib(U)$ is in the interior of the dual polygon
$p^*$. It follows that $i\sum_{i=1}^n W_{e_i}\geq 0$, with equality if and
only if the first-order deformation $U$ is trivial. The end of the proof is
the same as in the Euclidean case.
\medskip

Note that the rigidity argument given here for equivariant polyhedral
embeddings in $\R^3_1$ could perhaps be used in other situations, for instance
for closed, convex polyhedra in $\R^3_1$ which have faces which are not
necessarily space-like. However the infinitesimal rigidity of convex polyhedra
in $\R^3_1$ -- even when some faces are not space-like -- follows from the
infinitesimal rigidity of convex polyhedra in the Euclidean space (see
e.g. \cite{cpt}) -- so that having a direct Minkowski proof is not so
important. 

\section{Natural metrics on the moduli space of polygons}

The spaces of polygons in the sphere (or in other constant curvature spaces)
are of interest in topological or algebraic terms (see
e.g. \cite{kapovich-millson,kapovich-millson-sphere}) but also in metric
terms. This is in particular true of spaces of convex Euclidean polygons with
fixed angles (rather than fixed edge lengths); Bavard and Ghys
\cite{bavard-ghys} showed that those spaces can be naturally identified with
hyperbolic polyhedra, which are Coxeter polyhedra when some fairly simple
conditions on the angles of the polygons are satisfied. This is related to
the construction by Thurston \cite{SOP} of complex hyperbolic orbifolds as
moduli spaces of flat metrics with conical singularities on the sphere (see
\cite{fillastre1}), and also related to \cite{deligne-mostow}.

The positivity property of $b$ stated in Theorem B$_S$ can be used to define
some natural metrics on the spaces of convex polygons in the sphere (or in the
hyperbolic plane) with given edge lengths. We will describe those
constructions here, and indicate how they appear to be 
related to the metric used by
Bavard and Ghys on spaces of Euclidean polygons with given angles, which is
recovered as a limit case. The finer properties of the metrics on the spaces
of spherical or hyperbolic polygons are not studied here. 

\pg{Definitions of some metrics.}

The definitions of the metrics we want to consider stem from the following
elementary remark. We call $P^c_{n,S}$ the space of convex polygons with $n$
vertices in $S^2$, and $\cP^c_{n,S}$ the quotient of $P^c_{n,S}$ by $SO(3)$.
Given an $n$-uple $(l_1, \cdots, l_n)\in \R_+$, we
call $P^c_S(l)$ the space of convex polygons in $S^2$ with edge lengths
equal to $l$, and $\cP^c_S(l):=P^c_S(l)/SO(3)$.

\begin{remark} \label{rk:metrique}
Let $F:P^c_{n,S}\rightarrow S^2$ be a map such that:
\begin{itemize}
\item For all $p\in P^c_{n,S}$, $F(p)$ is in the interior of $p$.
\item For all $p\in P^c_{n,S}$ and $\gamma \in SO(3)$, $F(\gamma p)=\gamma
  F(p)$.
\end{itemize}
Let $p_0\in P^c_{n,S}$ and let $l$ be its edge lengths. 
Then the bilinear form $\langle b_2(\cdot, \cdot), F(p)\rangle$ on the tangent
space to $P^c_{S}(l)$ at each point $p$ defines a Riemannian metric on
$P^c_{S}(l)$, which is compatible with the quotient by $SO(3)$. Therefore it
defines a Riemannian metric on $\cP^c_{S}(l)$.
\end{remark}

\begin{proof}
By Theorem B$_S$, for each $p\in P^c_{S}(l)$ and each non-trivial
deformation $U\in
T_pP^c_{S}(l)$, $b(U)\in \R_+^* \mbox{int}(p^*)$. 
This means precisely that the scalar 
product of $b(U)$ with any point in the interior of $p$ is positive, and
this holds in particular of the scalar product of $b(U)$ with $F(p)$, so
that $\langle b_2(\cdot, \cdot), F(p)\rangle$ defines a positive semi-definite
bilinear form on $P^c_S(l)$, and the kernel corresponds precisely to the
trivial deformations.  

We have already seen that, if $U'$ is a trivial first-order deformation of
$p$, then $b(U+U')=b(U)$. The trivial first-order deformations of $p$, seen as
vectors in $T_pP^c_{S}(l)$, 
are the elements of the kernel of the differential
of the projection from $T_pP^c_{S}(l)$ to 
$T_{[p]}\cP^c_{S}(l)$, where
$[p]$ is the image of $p$ under the quotient of $P^c_{S}(l)$ by $SO(3)$. So
the Riemannian metric defined by $\langle b_2(\cdot, \cdot), F(p)\rangle$ on
$T_pP^c_{S}(l)$ is the pull-back by the projection of a Riemannian metric on
$T_{[p]}\cP^c_{S}(l)$.

Now let $\gamma\in SO(3)$, and let $U$ be a first-order deformation of
$p$. Considering $\gamma$ as a map acting on $P^c_S(l)$, we
associate to $U$ a first-order deformation $(d_p\gamma)(U)$ of $\gamma
p$. The definition of $b$ shows that it is ``invariant'' under the action of
$\gamma$, i.e. that $b((d_p\gamma)(U))=\gamma b(U)$. It follows that:
$$ \langle b((d_p\gamma)(U)), F(\gamma p)) \rangle = \langle \gamma b(U),
\gamma F(p)\rangle = \langle b(U), F(p)\rangle~, $$
so that the bilinear form $\langle b_2(\cdot, \cdot), F(p)\rangle$ is
invariant under the action of $SO(3)$.
\end{proof}

There are several possible choices for the function $F$. It is quite natural
to take some kind of barycenter of $p$, for instance:
\begin{itemize}
\item The barycenter of the vertices of $p$ (all with the same weight), which
  we will call $C_v(p)$. This
  is the most obvious solution, however we will see below that it lacks one
  desirable property. 
\item The barycenter of the area form on the 
  interior of $p$, which we call $C_i(p)$. 
  This is another quite obvious
  choice, we will see below that it has some better properties.
\item The barycenter of the vertices of $p$, with weights equal to the
  exterior angles, which we call $C_\alpha(p)$.
\item The barycenter of $p$ itself, i.e. of the union of its edges, 
  denoted here by $C_\partial(P)$.
\end{itemize}
Given a family of lengths $l=(l_1, \cdots, l_n)$, we will call $g_v(l)$
(resp. $g_i(l)$, $g_\alpha(l)$, $g_\partial(l)$) the Riemannian metrics on 
$\cP^c_S(l)$ described by Remark
\ref{rk:metrique} with $F$ equal to $C_v$ (resp. to $C_i$, $C_\alpha$,
$C_\partial$).

\pg{A geometric property of the barycenter of the area.}

The drawback with the choice of the barycenter of the vertices of $p$ is that
it does not have the following property, which is satisfied by the barycenter
of the interior of $p$.

\begin{prop} \label{pr:barycenter}
  Let $p$ be a convex polygon in $S^2$. The barycenter $C_i(p)$ of the
  interior of $p$ is contained in the 
  interior of the dual polygon $p^*$.
\end{prop}

\begin{proof}
By definition, $C_i(p)$ is in the interior of 
$p^*$ if and only if, for all point $x$ in the
interior of $p$, $\langle C_i(p), x\rangle \geq 0$.
This is equivalent to the fact that $C_i(p)$ is at distance at most $\pi/2$
from any point in $p$. Since the space of convex polygons with $n$ vertices is
connected, the result follows from the following assertion: if $p$ is
contained in the hemisphere centered at $C_i(p)$ and one point of $p$ is in
the boundary of this hemisphere, then $p$ is a polygon with 2
edges, which are both one half of a great circle, with one endpoint on $x$.

So we suppose that $p$ is contained in the hemisphere $H$ centered at $C_i(p)$,
with one point $x$ in the boundary, and that it has more than two edges. 
Then the connected component of $x$ in the intersection of $p$ with
$\partial H$ is either only one point $x$, or it an edge of $p$.
We first consider the first case.  

Consider the projective model of the hemisphere $H$, obtained by projecting
the points of $H$ radially to the plane in $\R^3$ tangent to $S^2$ at
$C_i(p)$. The image of $p$ in this model is a ``non-compact'' convex polygon
which we call $p_e$,
i.e. the boundary of a non-compact polygonal domain, with two parallel
infinite edges $e$ and $e'$, which is the image of the interior of $p$.

The projective model has the property that the symmetry in $H$ with respect to
the geodesic line $g$ at distance $\pi/2$ from $x$ (which contains $C_i(p)$)
acts like the symmetry in
$\R^2$ with respect to the line $g_e$ containing the origin and orthogonal to
$e$ and to $e'$. But a simple convexity argument shows that, if $y\in \R^2$ is
on the same side of $g_e$ as the infinite ends of $e$ and $e'$, and if $y$ is
not in the interior of $p_e$, then the image 
$y'$ of $y$ under the symmetry with
respect to $g_e$ is not in the interior of $p_e$; otherwise, some support line
of $p_e$ would go between $y$ and $y'$, and then it would have to ``cut'' the
infinite edges $e$ and $e'$. Moreover, $p_e$ is not invariant under the
symmetry with respect to $g_e$, unless it is an infinite strip, and then $p$
is a polygon with two edges, each of which is half of a great circle.

It follows that the same statement is true in $H$: if a point $z$ is on the
same side of $g$ as $x$, but is not contained in the interior of $p$, then the
image of $z$ under the symmetry with respect to $g$ is not contained in the
interior of $p$. Thus:
$$ \int_{\mbox{int}(p)} \langle x, \cdot\rangle da >0~, $$
and this contradicts the fact that $C_i(p)$, which is orthogonal to $x$, is
the barycenter of the interior of $p$. This shows the result when $x$ is a
vertex of $p$.

If the intersection of $p$ with $\partial H$ contains an edge of $p$, the same
argument can be used, taking as $x$ one point of this edge. The only
difference is that the two infinite edges of $p_e$ are not parallel; the key
point remains that, if $H_+$ is the hemisphere bounded by $g$ and containing
$x$ and $H_-$ is the other hemisphere bounded by $g$, then 
$H_-\cap \mbox{int}(p)$
is contained in the image under the symmetry with respect to $g$ of $H_+\cap
\mbox{int}(p)$. The proof can again be obtained using the projective model of
$H$. We leave the details to the reader.
\end{proof}

Note that this statement is not correct with $C_i(p)$ replaced with $C_v(p)$,
however it might hold with $C_i(p)$ replaced by $C_\alpha(p)$ or
$C_\partial(p)$. 

\pg{Euclidean polygons in the limit.}

Let $\alpha=(\alpha_1,\cdots, \alpha_n)$ be a set of positive numbers with sum
equal to $2\pi$. We call $\cP_E^*(\alpha)$ the space of convex Euclidean
polygons with exterior angles equal to $\alpha_1, \cdots, \alpha_n$,
considered up to the Euclidean isometries. We denote by $\cP^{*,A}_E(\alpha)$
the subspace of $\cP_E^*(\alpha)$ of polygons of area equal to $A$.

There is a natural metric on $\cP^{*,1}_E(\alpha)$, 
defined by Bavard and Ghys \cite{bavard-ghys}, 
called the ``area form'', which we will denote by $g_A$ here. 
It can be defined by considering the convex polygons
with edges parallel to some given directions, up to translation 
(this is very close to considering polygons with given angles, up to
isometry), and noting that the area is a quadratic form over the space of
those polygons. Bavard and Ghys show that the signature of this form is
$(1,n-3)$, and the moduli space of Euclidean polygons with given angles (up to
the homotheties), with the induced metric, is isometric to the interior of a
finite volume hyperbolic polyhedron. Moreover, under some explicit
conditions on $\alpha$, this polyhedron is Coxeter, i.e. the group generated
by the reflections in its faces is discrete. There is a relationship between
those metrics on spaces of convex Euclidean polygons with given angles and
the metrics defined above on spaces of spherical polygons with given edge
lengths. 

\begin{thm} \label{tm:convergence}
Let $(l^k)_{k\in \N}=(l_1^k, \cdots, l_n^k)_{k\in \N}$ 
be a sequence of $n$-uples
of positive numbers. 
Suppose that $\lim_{k\rightarrow \infty} l^k=\alpha$, with
$\alpha=(\alpha_1,\cdots, \alpha_n)$ such that $\sum_i \alpha_i=2\pi$
and that there is no $p,q\in \{
1,\cdots, n\}$ such that $\sum_{i=p}^q 
\alpha_i = \pi$. Then (by \cite{bavard-ghys}), $\cP^c_E(\alpha)$ is compact,
and: 
$$ (\cP^c_S(l^k), \frac{1}{2a_k} g_i(l^k))\rightarrow (\cP^{*,1}_E(\alpha), 
g_A)~, $$ 
where $a_k := 2\pi - \sum_{i=1}^n l_i^k$, and 
the convergence is Lipschitz. 
\end{thm}

The proof is given below, it uses some preliminary statements. 

\begin{remark} \label{rk:dual}
Let $l=(l_1,\cdots, l_n)$ be a family of positive numbers. 
\begin{enumerate}
\item Suppose that there
exists a convex polygon in $S^2$ with edge lengths given by $l$, then the sum
of the $l_i$ is less than $2\pi$.
\item If $p$ is a convex 
spherical polygon with edge lengths given by $l$, then $p^*$ has exterior
angles given by $l$, and area equal to $2\pi-\sum l_i$.
\item Suppose that $(p_k)_{k\in \N}$ is a sequence of convex spherical 
  polygons, with $p_k$ of
  edge lengths $l^k=(l_1^k,\cdots, l_n^k)$ with $\lim_{k\rightarrow \infty}
  l_i^k = l_i$ and $\sum_{i=1}^n
  l_i=2\pi$. Then, after taking a subsequence, $(p_k)_{k\in \N}$, considered
  as a sequence of subsets of $S^2$, converges
  either to a great circle, or to a polygon with two edges, each being one
  half of a great circle. This second case can happen only if there exist
  $p,q\in \{ 1,\cdots, n\}$ such that $\sum_{i=p}^q l_i=\pi$.
\end{enumerate}
\end{remark}

\begin{proof}
The proof of the first two points is either elementary or well-known. The
third point follows from the second, because the sequence of dual polygons
has area going to $0$, so that, after taking a subsequence, 
the dual polygons converge either to a point,
or to a segment. In the first case, $(p_k)$ converges to a great circle,
while, in the second case, $(p_k)$ converges to a polygon with two edges, both
of which is one half of the great circle which is dual to one of the endpoints
of the segment. 
If the convergence is to a segment, then the limits of the
exterior angles of the dual polygon 
at the vertices converging to each vertex of the segment has
to sum to $\pi$, and the condition on the lengths of the edges follows.
\end{proof}

We now need an additional notation concerning spaces of polygons with
fixed angles in the sphere. Given $\alpha=(\alpha_1, \cdots,
\alpha_n)$, we call $\cP_S^*(\alpha)$ the space of convex polygons in $S^2$
with angles given by $\alpha$, up to the global isometries. There is a natural
map $\Delta$ from $\cP_S^c(\alpha)$ to $\cP_S^*(\alpha)$ (and conversely),
defined by sending a convex polygon $p$ to the dual polygon $p^*$. Therefore,
we can consider on $\cP_S^*(\alpha)$ the pull-back metric $\Delta^* g_i$,
which we will call $g_i^*$.

It is actually helpful to consider $g_i$ and $g_i^*$ as symmetric bilinear
forms on the space of all convex polygons in $S^2$ with $n$ vertices,
$\cP^c_{S,n}$. Then $g_i$ is defined, as for the more restricted space of
polygons with given edge lengths, as:
$$ g_i(U,V) = \frac{1}{2} \sum_{i=1}^n \langle d\alpha_i(U)dv_i(V) +
d\alpha_i(V) dv_i(U), C_i(p)\rangle~, $$
and $g_i^*:=\Delta^*g_i$. In the same way, we define the area form, already
mentioned above, as a symmetric bilinear form over the space of convex 
Euclidean
polygons with $n$ vertices. For any point $x_0$ in the interior of $p$, the
area form can be defined as:
$$ g_A(U,V) = \frac{1}{4} \sum_{i=1}^n dl_i(U) dh_i(V) + dl_i(V) dh_i(U)~, $$
where $h_i$ is the distance from $x_0$ to the $i^{th}$ edge. 

Let $q=(v_1, \cdots, v_n)$ 
be a convex polygon in $S^2$. Let $u\in SO(3)$ be an isometry sending
its barycenter $C_v(q)$ to the ``north pole'' 
$(0,0,1)\in S^2$, and let $\rho(q)$ be the image of $u(q)$ by the projective
map from the upper hemisphere to $\R^2$. If $q$ is such that $C_v(q)$ is in
the interior of $q^*$, then  
$u(q)$ is contained in the upper hemisphere, so that $\rho(q)$ is a convex
Euclidean polygon. 
It is easy to check that $\rho(p)$, considered up
to the Euclidean isometries, is independent of the choice of $u$ (which is
uniquely defined up to a rotation).
We call $\rhob$ the map sending $q$ to $\rho(q)$. Moreover, given $q'\in
\cP^c_E$, we call $(\alpha_1^E, \cdots, \alpha_n^E)$ its angles, so that the
$\alpha_i^E$, for $1\leq i\leq n$, define $n$ functions on $\cP^c_E$.

\bprop \label{pr:cvbilin}
Let $(\alpha_1, \cdots, \alpha_n)$ be such that $\sum_i \alpha_i=2\pi$, but
such that there is no $p,q\in \{ 1,2, \cdots, n\}$ with $\sum_{i=p}^q \alpha_i
=\pi$. Let $(\alpha^k)_{k\in \N}$ be a sequence of 
$n$-uples of positive numbers
such that, for all $k\in \N$, $\sum_i\alpha_i^k<2\pi$,
and that $\lim_{k\rightarrow \infty} \alpha^k=\alpha$. 
Then there exists a sequence $(\epsilon_k)_{k\in
  \N}\rightarrow 0$ such that, for all $k\in \N$:
\begin{enumerate}
\item over $\cP^*_S(\alpha^k)$, 
$(1-\epsilon_k) g_i^* \leq 2\rhob^* g_A \leq (1+\epsilon_k) g_i^*$.
\item for all $1\leq i\leq n$, 
$\| d (\alpha^E_{i|\rhob(\cP^*_S(\alpha^k))})\|_{g_A}\leq \epsilon_k$.
\end{enumerate}
\eprop

\bpv
Let $q\in \cP^c_{S,n}$, and let 
$e_1, \cdots, e_n$ be the edges of $q$. Let $C_0:=C_i(q^*)$, we consider
here $C_0$ as a fixed point in $S^2$, even when considering first-order
deformations of $q^*$. For each $i\in \{1, 2,\cdots,
n\}$, let $\beta_i$ be the spherical distance from $C_0$ to
$e_i$. Consider a first-order deformation $U$ of $q$, it determines a
first-order deformation $U^*$ of the dual polygon $q^*$. If
$v_i^*$ is the vertex of $q^*$ which is dual to $e_i$, we have:
$$ \langle dv_i^*(U^*), C_i(q^*)\rangle = d(\sin(\beta_i))(U) = \cos(\beta_i)
d\beta_i(U)~. $$
Since the edge lengths of $q$ are the exterior angles of $q^*$, it follows
that, calling $l_i$ the edge lengths of $q$, we have:
$$ \langle b(U^*), C_i(q^*)\rangle = 
\sum_{i=1}^n \cos(\beta_i) dl_i(U) d\beta_i(U)~, $$
so that:
$$ g_i^*(U,U) = \sum_{i=1}^n \cos(\beta_i) dl_i(U) d\beta_i(U)~. $$

We now add to $U$ a trivial deformation, so as to obtain that the barycenter
$C_v(q)$ is fixed by $U$, and also that one of the vertices of $q$, for
instance $v_1$, is not moved. There exists a constant $C>0$ (depending on $n$)
such that, for all $1\leq i\leq n$:
\beq \label{eq:dv}
\| dv_i(U) \| \leq C\left(\sum_j |dl_j(U)|+|d\alpha_j(U)|\right)~. \eeq

Let $q^E:= \rho(q)$, let $l_i^E$ be the edge lengths of $q^E$, and let $h_i$
be the distances between the edges of $q^E$ and $\rho(C_0)$. Let $U^E$ be the
first-order deformation of $q^E$ which is induced by $U$ through $\rho$. 
We now suppose that
$q\in \cP^*_S(\alpha^k)$ (for some $k\in \N$). It follows from Remark
\ref{rk:dual} that, if $k$ is large enough, $q$ is contained in the ball
centered at $C_i(q^*)$ of radius $\lambda_k$, where $\lim_{k\rightarrow
  \infty} \lambda_k=0$. But the projective map $\rho$ is smooth, and its
differential at the ``north pole'' is an isometry. Thus
tt follows from (\ref{eq:dv}) that, for $k$
large enough, we have, for some constant $C>0$:
\beq \label{eq:dl} 
| dl^E_i(U^E) - dl_i(U) | \leq C  \lambda_k \left(\sum_j
  |dl_j(U)|+|d\alpha_j(U)|\right)~, \eeq
\beq \label{eq:dh} | dh_i(U^E) - d\beta_i(U) |\leq C \lambda_k \left(\sum_j
  |dl_j(U)|+|d\alpha_j(U)|\right)~. \eeq
Moreover, the same argument -- using the fact that $d\rho$ is conformal at the
``north pole'' -- yields that: 
\beq \label{eq:dalpha} 
| d\alpha^E_i(U^E) - d\alpha_i(U) | \leq C \lambda_k \left(\sum_j
  |dl_j(U)|+|d\alpha_j(U)|\right)~. \eeq
Using the description of $g_i^*$ and the definition of
$g_A$ given above, we obtain with (\ref{eq:dl}) and (\ref{eq:dh}) that, for
another constant $C'>0$: 
\beq \label{eq:gcomp}
|g_A(U^E,U^E) - g_i^*(U,U)| \leq C' \lambda_k \left(\sum_j
  |dl_j(U)|+|d\alpha_j(U)|\right)^2~. \eeq

On the other hand, we have seen in section 3 that $\langle b, C_i(q^*)\rangle$
is positive definite on the tangent plane at each point of
$\cP^c_S(\alpha^k)$; it follows by compactness that there exists a constant
$C''>0$ such that, for all $k\in \N$, all $q\in \cP^*_S(\alpha^k)$ and all
first-order deformation $U$ of $q$ leaving the angles constant at first order:
$$ g_i^*(U,U) \geq \frac{1}{C''} \left(\sum_j
  |dl_j(U)|+|d\alpha_j(U)|\right)^2~. $$
The first point of the proposition follows from this equation and from
(\ref{eq:gcomp}), while the second point then follows from (\ref{eq:dalpha}). 
\epv

\begin{proof}[Proof of Theorem \ref{tm:convergence}]
The first point of Proposition \ref{pr:cvbilin} shows that the restriction of
the $\rhob$ to $\cP^*_S(\alpha^k)$ is a Lipschitz map, which multiplies the
metric by $2$ (up to a factor converging to $1$). The polygons which are
in the image are, by Remark \ref{rk:dual}, close to $0$ in $\R^2$. Since
$\rho$ is an isometry at the ``north pole'', it follows that their area is
equivalent to $a_k$ as $k\rightarrow \infty$. 

Let $H_k$ be the homothety of ratio $1/\sqrt{a_k}$, acting on Euclidean
polygons. It follows from the above argument that $H_k\circ \rhob
(\cP^*_S(\alpha^k))$ converges to $\cP^{*,1}_E(\alpha)$. Moreover, the second
point of Proposition \ref{pr:cvbilin} shows that the norm of the differential
of the functions $\alpha_i^E$ converges to $0$ on the image, so that the
$H_k\circ \rhob (\cP^*_S(\alpha^k))$ actually converge $C^1$ to
$\cP^{*,1}_E(\alpha)$, and therefore the metric induced on the $H_k\circ \rhob
(\cP^*_S(\alpha^k))$ by $g_A$ converge to the metric on
$\cP^{*,1}_E(\alpha)$. 

But the first point of Proposition \ref{pr:cvbilin} shows that the metrics 
 induced on the $H_k\circ \rhob (\cP^*_S(\alpha^k))$ are close to $(1/2a_k)
 g^*_i$, and the proof of the theorem follows.
\end{proof}

A statement analog to Theorem \ref{tm:convergence}
appears to hold when there exist $p,q\in \{
1,\cdots, n\}$ such that $\sum_{i=p}^q 
\alpha_i = \pi$. In that case, however, $\cP^{*,1}_E(\alpha)$ 
has at least one ideal vertex, and the convergence should only be 
from the complement of some domains of 
$(\cP^c_H(l^k), \frac{k}{2} g_i)$ (corresponding to ``degenerate'' polygons
which are close to segments) to the complement of some neighborhoods of
the ideal vertices of $\cP^c_E(\alpha)$.

\pg{Hyperbolic and de Sitter polygons.}

Similar natural metrics can be constructed on the space of convex hyperbolic
polygons of given edge lengths, 
using Theorem B$_H$ rather than Theorem B$_S$. However, the
convergence to the spaces of Euclidean polygons should be considered for
spaces of hyperbolic polygons with fixed angles, rather than fixed edge
lengths. Fixing the angles of convex hyperbolic polygons is equivalent to
fixing the edge lengths of the dual polygons in the de Sitter plane (see
e.g. \cite{coxeter-de-sitter,coxeter-non-euclidean}). 

For convex polygons in the de Sitter plane, Theorem B$_{dS}$ (as seen 
in section 6) can also be used
to define natural metrics on the space of de Sitter polygons with given edge
lengths. Now, the sum of the lengths of the edges of a convex de Sitter
polygon (which is dual to a hyperbolic polygon)
is more than $2\pi$; the arguments used in the proof of Theorem 
\ref{tm:convergence} have analogs for the de Sitter plane, which seem to
indicate that, as the sum of the edge 
lengths converges to $2\pi$, the spaces of de
Sitter polygons, with adequately normalized
metrics, converge to a space of Euclidean
polygons, with the metric defined by \cite{bavard-ghys}.

\section*{Acknowledgements}

I would like to thank Philippe Eyssidieux, Vincent Guirardel,
Sergiu Moroianu and Boris Springborn for some helpful discussions and
interesting comments related to the content of this paper. 

\bibliographystyle{alpha}

\begin{thebibliography}{Thu98}

\bibitem[Ale58]{alex}
A.~D. Alexandrow.
\newblock {\em Konvexe polyeder}.
\newblock Akademie-Verlag, Berlin, 1958.

\bibitem[BG92]{bavard-ghys}
Christophe Bavard and {\'E}tienne Ghys.
\newblock Polygones du plan et poly\`edres hyperboliques.
\newblock {\em Geom. Dedicata}, 43(2):207--224, 1992.

\bibitem[Cau13]{cauchy}
A.~L. Cauchy.
\newblock Sur les polygones et poly\`edres, second m\'emoire.
\newblock {\em Journal de l'Ecole Polytechnique}, 19:87--98, 1813.

\bibitem[Cox43]{coxeter-de-sitter}
H.~S.~M. Coxeter.
\newblock A geometrical background for de {S}itter's world.
\newblock {\em Amer. Math. Monthly}, 50:217--228, 1943.

\bibitem[Cox57]{coxeter-non-euclidean}
H.~S.~M. Coxeter.
\newblock {\em Non-{E}uclidean geometry}.
\newblock Mathematical Expositions, no. 2. University of Toronto Press,
  Toronto, Ont., 1957.
\newblock 3rd ed.

\bibitem[Cox93]{coxeter-projective}
H.~S.~M. Coxeter.
\newblock {\em The real projective plane}.
\newblock Springer-Verlag, New York, third edition, 1993.
\newblock With an appendix by George Beck, With 1 IBM-PC floppy disk (5.25
  inch; DD).

\bibitem[Deh16]{dehn-konvexer}
M.~Dehn.
\newblock {\"Uber} den {Starrheit} konvexer {Polyeder}.
\newblock {\em Math. Ann.}, 77:466--473, 1916.

\bibitem[DM86]{deligne-mostow}
P.~Deligne and G.~D. Mostow.
\newblock Monodromy of hypergeometric functions and nonlattice integral
  monodromy.
\newblock {\em Inst. Hautes \'Etudes Sci. Publ. Math.}, (63):5--89, 1986.

\bibitem[Euc02]{euclid}
Euclid.
\newblock {\em Elements}.
\newblock Green Lion Press, Santa Fe, NM, 2002.
\newblock All thirteen books complete in one volume, The Thomas L. Heath
  translation, Edited by Dana Densmore.

\bibitem[Fil]{fillastre1}
F.~Fillastre.
\newblock Spaces of polygons to spaces of polyhedra following {Bavard}, {Ghys}
  and {Thurston}.
\newblock math.MG/0308187.

\bibitem[Fil92]{filliman}
P.~Filliman.
\newblock Rigidity and the {A}lexandrov-{F}enchel inequality.
\newblock {\em Monatsh. Math.}, 113(1):1--22, 1992.

\bibitem[Isk00]{iskhakov}
I.~Iskhakov.
\newblock {\em On hyperbolic surface tessellations and equivariant spacelike
  convex polyhedral surfaces in Minkowski space}.
\newblock PhD thesis, Ohio State University, 2000.

\bibitem[Kan90]{kann}
Edgar Kann.
\newblock Infinitesimal rigidity of almost-convex oriented polyhedra of
  arbitrary {E}uler characteristic.
\newblock {\em Pacific J. Math.}, 144(1):71--103, 1990.

\bibitem[KM95]{kapovich-millson}
Michael Kapovich and John Millson.
\newblock On the moduli space of polygons in the {E}uclidean plane.
\newblock {\em J. Differential Geom.}, 42(2):430--464, 1995.

\bibitem[KM99]{kapovich-millson-sphere}
Michael Kapovich and John~J. Millson.
\newblock On the moduli space of a spherical polygonal linkage.
\newblock {\em Canad. Math. Bull.}, 42(3):307--320, 1999.

\bibitem[LegII]{legendre}
A.-M. Legendre.
\newblock {\em El\'ements de g\'eom\'etrie}.
\newblock Paris, 1793 (an II).
\newblock Premi\`ere \'edition, note XII, pp.321-334.

\bibitem[LS00]{iie}
Fran\c{c}ois Labourie and Jean-Marc Schlenker.
\newblock Surfaces convexes fuchsiennes dans les espaces lorentziens \`a
  courbure constante.
\newblock {\em Math. Annalen}, 316:465--483, 2000.

\bibitem[Mil71]{milka-polygones}
A.~D. Milka.
\newblock The isoperimetry of polygons on the sphere.
\newblock {\em Ukrain. Geometr. Sb.}, (10):49--50, 1971.

\bibitem[O'N83]{O}
B.~O'Neill.
\newblock {\em Semi-Riemannian Geometry}.
\newblock Academic Press, 1983.

\bibitem[Pog56]{pogo-polygones}
A.~V. Pogorelov.
\newblock A new proof of rigidity of convex polyhedra.
\newblock {\em Uspehi Mat. Nauk (N.S.)}, 11(5(71)):207--208, 1956.

\bibitem[RR00]{rodriguez-rosenberg}
Lucio Rodr{\'{\i}}guez and Harold Rosenberg.
\newblock Rigidity of certain polyhedra in {${\bf R}\sp 3$}.
\newblock {\em Comment. Math. Helv.}, 75(3):478--503, 2000.

\bibitem[Sch98]{shu}
Jean-Marc Schlenker.
\newblock M\'etriques sur les poly\`edres hyperboliques con\-vexes.
\newblock {\em J. Differential Geom.}, 48(2):323--405, 1998.

\bibitem[Sch01]{cpt}
Jean-Marc Schlenker.
\newblock Convex polyhedra in {Lorentzian} space-forms.
\newblock {\em Asian J. of Math.}, 5:327--364, 2001.

\bibitem[Sie02]{siegel}
A.~Siegel.
\newblock A {D}ido problem as modernized by {F}ejes {T}\'oth.
\newblock {\em Discrete Comput. Geom.}, 27(2):227--238, 2002.

\bibitem[Ste42]{steiner}
J.~Steiner.
\newblock Sur le maximum et le minimum des figures dans le plan, sur la
  sph\`ere, et dans l'espace en g\'en\'eral.
\newblock {\em J. Reine Angew. Math.}, 24:93--152, 190--250, 1842.

\bibitem[Sto68]{stoker}
J.~J. Stoker.
\newblock Geometrical problems concerning polyhedra in the large.
\newblock {\em Comm. Pure Appl. Math.}, 21:119--168, 1968.

\bibitem[Thu98]{SOP}
William~P. Thurston.
\newblock Shapes of polyhedra and triangulations of the sphere.
\newblock In {\em The Epstein birthday schrift}, volume~1 of {\em Geom. Topol.
  Monogr.}, pages 511--549 (electronic). Geom. Topol. Publ., Coventry, 1998.

\bibitem[Vol56]{volkov-polygones}
Yu.~A. Volkov.
\newblock On deformations of a convex polyhedral angle.
\newblock {\em Uspehi Mat. Nauk (N.S.)}, 11(5(71)):209--210, 1956.

\end{thebibliography}

\end{document}